\documentstyle[amscd]{amsart}

\def\AA{{\blb A}}
\def\CC{{\blb C}}

\def\LL{{\blb L}}
\def\NN{{\blb N}}

\def\PP{{\blb P}}

\def\ZZ{{\blb Z}}

\def\blb#1{\Bbb#1}

\def\11{{1\kern-3.5pt 1}}
\def\mumu{{\mu\kern-4.2pt\mu}}

\def\boxtimes{\setbox0\hbox{$\Box$}\copy0\kern-\wd0\hbox{$\times$}}

\def\Uol{{\bar U}}
\def\Vol{{\bar V}}

\def\End{\operatorname {End}}

\def\gr{\operatorname {gr}}

\def\Hom{\operatorname {Hom}}

\def\ker{\operatorname {ker}}

\def\op{{\operatorname {op}}}

\def\Spec{\operatorname {Spec}}

\def\th{\operatorname {th}}    

\def\Ann{\operatorname{Ann}}

\def\Comod{{\sf Comod}}

\def\dim{\operatorname{dim}}

\def\End{\operatorname{End}}

\def\Ext{\operatorname{Ext}}

\def\Fdim{{\sf Fdim}}

\def\fl.dim{\operatorname{flat.dim}}

\def\GrMod{{\sf GrMod}}

\def\Hom{\operatorname{Hom}}

\def\id{\operatorname{id}}

\def\Kdim{\operatorname{Kdim}}

\def\Lex{{\sf Lex}}

\def\max{\operatorname{max}}

\def\mod{{\sf mod}}
\def\Mod{{\sf Mod}}

\def\op{\operatorname{op}}

\def\Proj{{\sf Proj}}

\def\smallbullet{\scriptscriptstyle{\bullet}}

\def\Spec{\operatorname{Spec}}
\def\Sum{{\sf Sum}}

\def\Tails{{\sf Tails}}

\def\uExt{\operatorname{\underline{Ext}}}
\def\cExt{\operatorname{\cal{E}{\it xt}}}

\def\l{\leftarrow}

\def\d{\downarrow}

\let\oldtext\text
\def\text#1{\oldtext{\normalshape #1}}

\def\a{\alpha}
\def\b{\beta}
\def\c{\gamma}
\def\d{\delta}

\def\ve{\varepsilon}

\def\l{\lambda}
\def\s{\sigma}

\def\fm{{\frak m}}

\def\fsl{{\frak sl}}

\def\sA{{\sf A}}

\def\sS{{\sf S}}
\def\sT{{\sf T}}

\def\coh{{\sf{ coh}}}

\def\Qcoh{{\sf{ Qcoh}}}

\def\cB{{\cal B}}

\def\cD{{\cal D}}

\def\cF{{\cal F}}
\def\cG{{\cal G}}

\def\cHom{{\cal H}{\it om}}
\def\cI{{\cal I}}
\def\cJ{{\cal J}}

\def\cL{{\cal L}}
\def\cM{{\cal M}}
\def\cN{{\cal N}}
\def\cO{{\cal O}}

\def\cS{{\cal S}}

\def\cTor{{\cal T}\!{\it or}}

\def\cstar{\operatorname {!}}

\newtheorem{lemma}{Lemma}[section]
\newtheorem{proposition}[lemma]{Proposition}
\newtheorem{theorem}[lemma]{Theorem}
\newtheorem{corollary}[lemma]{Corollary}

{

}

\theoremstyle{definition}

\newtheorem{example}[lemma]{Example}
\newtheorem{definition}[lemma]{\sl Definition}

\theoremstyle{remark}

\numberwithin{equation}{section}

\begin{document}

\pagenumbering{arabic}

\title{Subspaces of non-commutative spaces}

\author{ S. Paul Smith}

\address{Department of Mathematics, Box 354350, Univ.
Washington, Seattle, WA 98195, USA}

\email{smith@@math.washington.edu}

\thanks{The author was supported by NSF grant DMS-0070560 }

\keywords{}

\begin{abstract}
This paper concerns the closed points, closed subspaces, open
subspaces, weakly closed and weakly open
subspaces, and effective divisors, on a non-commutative space.
\end{abstract}

\maketitle

\section{Introduction}

This paper examines non-commutative 
analogues of some of the elementary ideas of topology and 
(algebraic) geometry.
Our wish is to have a language for non-commutative
algebraic geometry that is {\it geometric} and {\it topological}.
This will, eventually, allow us to use our geometric intuition in
situations that at present are viewed almost exclusively from an
algebraic perspective.
One can already see hints of this goal in Gabriel's thesis
\cite{G}. More recently, Rosenberg and Van den Bergh 
have given new life to this program. 
We further extend their ideas by examining 
the interactions between the following basic objects: weakly closed 
and weakly open subspaces, closed subspaces, open subspaces, 
closed points, and effective divisors.

We follow Rosenberg \cite{Rosen} \cite{Ros2} and Van den Bergh 
\cite{vdB} in taking a Grothendieck
category as our basic non-commutative geometric object.  Thus a
non-commutative space $X$ is a Grothendieck category $\Mod X$. 
The idea for the notation $X=\Mod X$ is Van den Bergh's; we also
use his notion of $X$-$X$-bimodules, and denote the unit object in
this monoidal category by $o_X$.

The standard commutative example of a space
is the category of quasi-coherent sheaves on a quasi-separated,
quasi-compact scheme. 
If $X$ is such a scheme we will speak of it as a space with the
understanding that $\Mod X$ is $\Qcoh X$.
The two non-commutative models are $\Mod R$, the category of right
modules over a ring, and $\Proj A$, the non-commutative projective
spaces having a (not necessarily commutative) graded ring $A$ as
homogeneous coordinate ring
(Definition \ref{defn.proj.space}).

The notion of a weakly closed subspace of a non-commutative space
first appeared in Gabriel's thesis 
\cite[page 395]{G}---he calls a full subcategory of an abelian
category closed if it is closed under subquotients
and direct limits. We will modify his terminology,
and say that a full subcategory of a Grothendieck category is
weakly closed if it satisfies these conditions. Closure under
subquotients ensures that the inclusion functor is exact, and 
closure under direct limits ensures that it has a right
adjoint. If the inclusion functor also has a left adjoint we say 
that the subcategory is closed. This modification of Gabriel's
terminology makes it compatible with the language of
algebraic geometry: the closed subcategories (in our sense)
of $\Mod R$ are in bijection with the two-sided ideals of $R$; and,
every closed subscheme of a scheme $X$ determines a
closed subspace of $\Qcoh X$. 
If $X$ is a quasi-projective scheme over a field then the closed
subspaces are the same things as the closed 
subschemes (Theorem \ref{thm.closed=closed}).
Our terminology is
also compatible with Van den Bergh's notion of weak ideals; weakly
closed subspaces of a space $X$ correspond to weak ideals in
$o_X$, the identity functor on $\Mod X$.

Proposition \ref{prop.bijection} shows that for noetherian
non-commutative spaces that are close to being commutative there is
a bijection between what we call the prime subspaces of $X$ and the
indecomposable injectives in $\Mod X$. This is extends Matlis'
result for commutative notherian rings, and Gabriel and 
Cauchon's extension of Matlis' result to rings 
satisfying Gabriel's condition (H) (fully bounded right 
noetherian rings).

We define the complement to a weakly closed subspace in the obvious
way (Definition \ref{defn.complement}), and define weakly 
open subspaces (Definition \ref{defn.open}) in such a way that they are
precisely the complements to the weakly closed subspaces when the
ambient space is locally noetherian.
The intersection of an arbitrary collection, and the union of a
finite collection, of weakly closed subspaces 
(Definitions \ref{defn.intersect}
and \ref{defn.union}) have already been defined by Rosenberg; these
are again weakly closed subspaces. Using the fact that a weakly open
subspace is a complement to a weakly closed subspace allows us
to define the intersection and union of weakly open subspaces. 
For the most part these
ideas interact with each other in reasonable ways---for example, 
see Proposition \ref{prop.W.cap.Z}, 
Lemmas \ref{lem.union.open}, \ref{lem.dist},
\ref{lem.UsubZ}, and \ref{lem.weak.closure}, Corollary
\ref{cor.open.in.open}, and Propositions \ref{prop.contains},
\ref{prop.Z.cap.U}, and \ref{prop.induced.top2}.
Furthermore, Proposition \ref{prop.fiber} shows
that if $U$ and $V$ are weakly open subspaces of a locally
noetherian space $X$,
then $U \cap V$ is the fiber product $U \times_X V$ in the category of
non-commutative spaces with maps being the morphisms.
Despite these positive results we do not know whether 
$W \cap(Y \cup Z)$ is equal to $(W \cap Y)
\cup (W \cap Z)$ when $W$, $Y$, and $Z$, are weakly closed subspaces;
the analogous formula for weakly open subspaces is valid.
Another problem is that if
$W$ and $Z$ are weakly closed subspaces
of $X$ such that $W \cap Z=\phi$, then $Z$ need not be contained
in $X \backslash W$; we determine exactly when $Z$ is
a weakly closed subspace of $X\backslash W$.
The fact that $Z$ need not be contained in $X \backslash W$ is an
essential feature of the non-commutativity of $X$. It reflects the
fact that there can be non-split extensions between non-isomorphic
simple modules.

In section \ref{sect.pts} we define closed points of a space.
A test case for our definition is the following. Let $A$ be a
connected graded $k$-algebra. Artin-Tate-Van den Bergh \cite{ATV2}
have defined the notion of a point module for $A$ with the idea
behind their definition being that point modules for $A$ should
play the role of points on the non-commutative projective space
having $A$ as a homogeneous coordinate ring. However, to date there
has been no definition of ``point'' in non-commutative geometry.
In section \ref{sect.pts} we show that point modules over $A$ give 
closed points of the associated projective space $\Proj A$.

To motivate the definition recall that a closed point $p$ on a scheme
$X$ corresponds to a simple object $\cO_p=\cO_X/\fm_p$ in $\Qcoh
X$, where
$\fm_p$ is the sheaf of ideals vanishing at $p$. This sets up a
bijection between closed points and simple objects in $\Qcoh X$.
Since the geometric objects in non-commutative geometry are
Grothendieck categories we define a closed point of a space $X$ to
be a certain kind of full subcategory of $\Mod X$. That subcategory 
is required to be closed, in the sense described above, and to
consist of all direct sums of a simple $X$-module, say $\cO_p$.
The requirement that the subcategory be closed means that
not all simple $X$-modules correspond to closed points. In
particular, an infinite dimensional simple module over a finitely
generated $k$-algebra does not correspond to a closed point; such
simple modules behave like higher dimensional geometric
objects (see the remarks after Proposition \ref{prop.proj.pts}).

Sometimes even finite dimensional simple modules behave 
like higher dimensional geometric objects.
The $\cD$-module approach to the representation
theory of semisimple Lie algebras, and the dimension of the
characteristic variety, yields a dimension function that gives finite
dimensional representations dimension greater than zero, thus
reflecting the fact that such modules do not always behave like
points. The behavior of the higher Ext-groups of finite dimensional 
simples also indicates that they behave like
higher dimensional geometric objects.

A further example of this phenomenon occurs when 
$\a:\tilde X \to X$ is Van den Bergh's blowup for a non-commutative 
surface $X$. Sometimes $X$ can have a
module of Krull dimension one whose strict
transform is a simple $\tilde X$-module, say $S$.
That simple module gives a closed point on $\tilde X$, but
with regard to its intersection properties (defined in terms of
$\Ext$-groups) that simple $\tilde X$-module behaves more like 
a curve than a point (see \cite{JS2}). 
This is reinforced by the fact that $\a_*S$ is the original 
module of Krull dimension one---in particular, $\a$ does not send
this closed point of $\tilde X$ to a closed point of $X$.  
For many purposes it is better to think
of $S$ as a geometric object of dimension one.

In Section \ref{sect.div} we examine effective divisors (as defined
by Van den Bergh) on a non-commutative space. These are analogues
of effective Cartier divisors; an effective divisor is defined 
in terms of an invertible ideal, say $o_X(-Y)$, in $o_X$. That ideal
cuts out a closed subspace $Y$ on $X$.
If $W$ is a weakly closed subspace of $X$, we give precise
conditions under which $W \cap Y$ is an effective divisor on $W$.
We also show that if $Y$ is ample, then $W \cap Y \ne \phi$ for all
$W$ of positive dimension. We give conditions under which the set
of closed points in $X$ is the disjoint union of the sets of closed
points in $Y$ and $X \backslash Y$.
Effective divisors behave better than arbitrary closed 
(or weakly closed) subspaces with respect to intersection and
containment of other subspaces.

\smallskip

{\bf Acknowledgements.}
We thank M. Artin, P. J\o rgensen, S. Kovacs, M. Van den Bergh, 
and J. Zhang for helpful conversations.
We thank C. Pappacena for asking questions that
prompted Propositions \ref{prop.fbn} and \ref{prop.bijection}.
This work was begun while the author was visiting the 
University of Copenhagen during the summers of 1998 and 1999, 
and was continued during
the workshop on Noncommutative Algebra at the Mathematical 
Sciences Research Institute at Berkeley in 2000. 
The author is grateful to both institutions for
their hospitality and financial support.

\section{Preliminaries}
\label{sect.qsch}

Throughout we work over a fixed commutative base ring $k$. 
All categories are $k$-linear and functors between them are $k$-linear.

\begin{definition}
A {\sf non-commutative space} $X$ is a Grothendieck category $\Mod X$.
Objects in $\Mod X$ are called $X$-modules.
We say $X$ is {\sf locally noetherian} if $\Mod X$ is locally 
noetherian (that is, if it has a set of noetherian generators).
\end{definition}

The {\sf empty space} $\phi$ is defined by declaring $\Mod \phi$ to
be the zero category; that is, the abelian category having only 
one object and one morphism.

We write $\Mod R$ for the category of right modules over a ring
$R$.

We say $X$ is {\sf affine} if $\Mod X$ has a progenerator, 
and in this case any ring $R$ for which $\Mod X$ is equivalent 
to $\Mod R$ is called a {\sf coordinate ring} of $X$.
It is easy to see that an affine space is locally noetherian if and 
only if one, and hence all, of its coordinate rings is right noetherian.

If $(X,\cO_X)$ is a scheme then the category $\Mod \cO_X$ 
of all sheaves of $\cO_X$-modules is a Grothendieck category. 
If $X$ is quasi-compact and quasi-separated (for example, if $X$ is
a noetherian scheme) 
the full subcategory of $\Mod \cO_X$ consisting of the 
quasi-coherent $\cO_X$-modules is a Grothendieck category
\cite[page 186]{G.SGA6}.
We denote this category by $\Qcoh X$. Whenever $X$ is a
quasi-compact and quasi-separated scheme we will speak of it as a
space in our sense with the tacit understanding that 
$\Mod X$ is $\Qcoh X$. For example, $\Spec \ZZ$ will denote the
space whose module category is $\Mod \ZZ$.

\medskip

The non-commutative analogues of projective schemes were first
studied in \cite{AZ} and \cite{Ver}.

\begin{definition}
\label{defn.proj.space}
Let $k$ be a field.
Let $A$ be an $\NN$-graded $k$-algebra such that $\dim_k A_n <
\infty$ for all $n$. Define $\GrMod A$ to be the category
of $\ZZ$-graded $A$-modules with morphisms the $A$-module homomorphisms
of degree zero. We write $\Fdim A$ for the full subcategory of direct limits 
of finite dimensional modules. We define the quotient category
$$
\Tails A=\GrMod A/\Fdim A,
$$
and denote by $\pi$ and $\omega$ the quotient functor and its right adjoint.
The {\sf projective space} $X$ with {\sf homogeneous coordinate ring}
$A$ is defined by $\Mod X:=\Tails A$. 

We define its structure module, $\cO_X$, to be $\pi A$, and we define 
$$
\Proj A=(X,\cO_X).
$$
This is an {\sf enriched} space in the sense of \cite[Section
3.6]{vdB}.

The {\sf cohomology groups} of a module $\cF$ over the enriched space
$X$ are defined in \cite{Ver} and \cite{AZ} to be
$$
H^q(X,\cF):=\Ext^q_X(\cO_X,\cF).
$$
\end{definition}

\medskip

The idea for the next definition is due to Rosenberg \cite[Section
1]{Ros2}.

\begin{definition}
\label{defn.map}
If $X$ and $Z$ are non-commutative spaces, 
a {\sf weak map} $f:Z \to X$
is a natural equivalence class of left exact 
functors $f_*:\Mod Z \to \Mod X$.
A weak map $f:Z \to X$ is a {\sf map} if $f_*$ has a left 
adjoint.
A weak map is {\sf affine} \cite[Section 4.3]{Ros2} 
if $f_*$ is faithful, and has a both
a left and a right adjoint.
A left adjoint to $f_*$ will be denoted by $f^*$, and a right
adjoint will be denoted by $f^!$. 
\end{definition}

A weak map $f$ is a map if and only if $f_*$ commutes
with products.

Let $f:X \to Y$ be a morphism of schemes. Suppose that $\Qcoh X$
and $\Qcoh Y$ are Grothendieck categories, so that $X$ and $Y$ are
spaces in our sense. If either $X$ is noetherian, or $f$ is 
quasi-compact and separated, then $f$ induces a map in the sense of
Definition \ref{defn.map} \cite[Proposition 5.8, Ch. II]{H}.
However, there are maps that are not induced by 
morphisms. For example, the inclusion
$$
\begin{pmatrix} k & 0 \\ 0 & k \end{pmatrix}
\to 
\begin{pmatrix} k & k \\ k & k \end{pmatrix}
$$
induces an affine map $\Spec k \to \Spec (k \times k)$ of spaces
that is not induced by a morphism between these schemes.

\begin{definition}
\label{defn.weakly.closed}
Let $X$ be a non-commutative space. 
A {\sf weakly closed subspace},
say $W$, of $X$ is a full subcategory $\Mod W$ of $\Mod X$ 
that is closed under subquotients (and therefore under isomorphisms), 
and for which the inclusion functor $\a_*:\Mod W \to \Mod X$ has 
a right adjoint, denoted by $\a^!$.
We write $\a:W \to X$ for the weak map corresponding to $\a_*$. 
A weakly closed subspace $W$ is {\sf closed} if
$\a_*$ also has a left adjoint, which is denoted by $\a^*$.
\end{definition}

For example, $X$ and $\phi$ are closed subspaces of $X$.

The inclusion of a closed subspace is an affine map.

Let $W$ be a weakly closed subspace of $X$ and let $\a:W \to X$ be
the inclusion. Then $\Mod W$ is a Grothendieck category and is
locally noetherian if $\Mod X$ is. 
Because $\Mod W$ is closed under subquotients,
$\a_*$ is an exact functor.
Because $\a^!$ has an exact left adjoint it sends injective
$X$-modules to injective $W$-modules.
The hypotheses imply that $\a^!\a_* \cong \id_W$ when $W$ is a
weakly closed, and $\a^*\a_* \cong \id_W$ when $W$ is closed.

When $Z$ is a closed subscheme of a scheme $X$ one often confuses
the structure sheaf $\cO_Z$ which lies in $\Qcoh Z$ with its direct 
image sheaf which lies in $\Qcoh X$.
We will do something similar when $W$ is a weakly closed subspace of a
space $X$. If $\a:W \to X$ is the inclusion, we will confuse the 
trivial $W$-$W$-bimodule $o_W$ with the weak $X$-$X$-bimodule 
corresponding to the left exact functor $\a_*\a^!$. Thus, we will 
speak of the map $o_X \to o_W$ of
weak $X$-$X$-bimodules, the kernel of which is the weak ideal
defining $W$.

\begin{definition}
\label{defn.open}
Let $X$ be a non-commutative space. 
A {\sf weakly open subspace}, say $U$,
of $X$ is a full subcategory $\Mod U$ of $\Mod X$ that is closed
under isomorphism and is such that the
inclusion functor $j_*:\Mod U \to \Mod X$ has an exact left adjoint
$j^*$.
If $V$ is another weakly open subspace of $X$, we say that $U$ is
{\sf contained in} $V$, and write $U \subset V$, if $\Mod U \subset
\Mod V$.
\end{definition}

If a weakly closed subspace is also a weakly open subspace, then
it is a closed subspace.

\begin{definition}
We call a weak map $f:Y \to X$ a {\sf closed immersion} 
(resp., {\sf weakly closed immersion}, {\sf weakly open immersion}, 
{\sf open immersion}) if it is an isomorphism onto a closed (resp.
weakly open, weakly closed, open) subspace of $X$.
\end{definition}

The next two results are standard: they follow from the fact that
$\a^!$ and $j_*$ are right adjoints to exact functors and so
preserve injectives.

\begin{lemma}
Let $\a:W \to X$ be a weakly closed immersion. If $M \in \Mod W$ and 
$N \in \Mod X$, then there is a spectral sequence
\begin{equation}
\label{eq.sp.seq}
\Ext^p_W(M,R^q\a^!N) \Rightarrow \Ext_X^{p+q}(\a_*M,N).
\end{equation}
\end{lemma}

\begin{lemma}
If $j:U \to X$ is a weakly open immersion, then
there is a spectral sequence
\begin{equation}
\label{eq.spec.seq.open}
\Ext^p_X(M,R^qj_*N) \Rightarrow \Ext^{p+q}_U(j^*M,N)
\end{equation}
whenever $M \in \Mod X$ and $N \in \Mod U$.
\end{lemma}

\section{Weakly closed subspaces}
\label{sect.wclosed}

Rosenberg has shown that the closed subspaces of an affine space
are in bijection with the two-sided ideals of its coordinate ring.
A typical non-commutative ring has few two-sided ideals, and a
typical non-commutative space has few closed subspaces.
A more useful notion for non-commutative geometry is that of a
weakly closed subspace. 

This section shows that there are reasonable notions of 
intersection and union for weakly closed subspaces.

The notion of a weakly closed subspace makes sense for
schemes. For example, $\Spec \ZZ$ has a weakly closed subspace 
$W$ defined by declaring $\Mod W$ to consist of the 
torsion abelian groups. There is a proper descending chain of
weakly closed subspaces $W=W_0 \supset W_1 \supset W_2 \supset
\ldots$ defined by declaring $\Mod W_i$ to be all torsion groups
that have no torsion with respect to the first $i$ positive primes.

If $X$ is a scheme such that $\Qcoh X$ is a Grothendieck category,
and $Z$ is a closed subscheme of $X$, then the
full subcategory of $\Qcoh X$ consisting of the modules whose
support lies in $Z$ is closed under subquotients 
and direct limits so is the module category of a weakly closed 
subspace of $X$; this subspace is rarely closed.

Let $X$ be a space and $M$ an $X$-module.
Let $W$ be the smallest weakly closed
subspace such that $M \in \Mod W$. Then $\Mod W$ consists of all
subquotients of all direct sums of copies of $M$. This subcategory
is closed under subquotients, so the inclusion functor is
exact; it is closed under direct
sums because direct sums are exact in $\Mod X$; it is therefore
closed under direct limits, so the inclusion functor has a right
adjoint.

\medskip

The following is a triviality.

\begin{lemma}
\label{lem.sub.sub}
If $Z$ is a weakly closed subspace of $X$ and $W$ is a weakly
closed subspace of $Z$, then $W$ is a weakly closed subspace of
$X$.
\end{lemma}

\begin{definition}
\label{defn.intersect}
\label{defn.contain}
Let $Z$ and $W$ be weakly closed subspaces of $X$.
We say that $Z$ {\sf lies on} $W$ or is 
{\sf contained in} $W$, if $\Mod Z$ is contained in $\Mod W$.
In this case we write $Z \subset W$.
Their {\sf intersection} $W \cap Z$ is defined by 
declaring that 
$$
\Mod W \cap Z := \Mod W \cap \Mod Z.
$$
That is, $\Mod W \cap Z$ is the full subcategory of $\Mod X$
consisting of those $X$ modules that belong to both $\Mod W$ and
$\Mod Z$.
\end{definition}

\begin{proposition}
\label{prop.W.cap.Z}
Suppose that $W$ and $Z$ are weakly closed subspaces of $X$.
Then
\begin{enumerate}
\item{}
$W$ and $Z$ are spaces;
\item{}
$W \cap Z$ is weakly closed in $X$, and in $W$, and in $Z$;
\item{}
if $Z$ is closed in $X$, then $W \cap Z$ is closed in $W$;
\item{}
if $W$ and $Z$ are closed in $X$, then $W \cap Z$ is closed in $X$,
and in $W$, and in $Z$;
\item{}
if $Z \subset W$, then $Z$ is a weakly closed subspace of $W$;
\item{}
if $Z \subset W$ and $Z$ is closed in $X$, then $Z$ is closed in
$W$;
\item{}
if $Z$ is also weakly open in $X$, then $W \cap Z$ is both closed
and weakly open in $W$.
\end{enumerate}
\end{proposition}
\begin{pf}
Let $\a:W \to X$ and $\b:Z \to X$ be the inclusions.
Let $\c_*:\Mod W \cap Z \to \Mod Z$ and  $\d_*:\Mod W \cap Z \to
\Mod W$ be the inclusion functors. 

(1)
To see that $W$ is a space one checks that $\Mod W$ is a Grothendieck 
category. By, for example, \cite[Proposition 3.4.3]{vdB}, $\Mod W$ 
is a cocomplete abelian category, and the inclusion $\a_*$ is exact 
and commutes with direct limits. Since direct limits in $\Mod W$
coincide with direct limits in $\Mod X$, $\Mod W$ has exact direct
limits. 

Since $\Mod X$ has a small set of generators,
the collection of all submodules of any $X$-module is a small set, and
hence so is the collection of all its subquotients. If
$\{M_i \; | \; i \in I\}$ is a small set of generators for $\Mod
X$, then the set of all subquotients of all $M_i$ that belong to
$\Mod W$ is small. This is a set of generators for $\Mod W$:
if $f:L \to N$ is a non-zero morphism in $\Mod W$, then there is a
morphism $g:M_i\to L$ such that $fg \ne 0$; since $L \in \Mod W$
$g$ induces a morphism $g':M' \to L$ from some subquotient 
$M'$ of $M_i$ belonging to $\Mod W$ such that $fg' \ne 0$.

(2)
Since $\Mod W$ and $\Mod Z$ are both closed under subquotients and
direct limits in $\Mod X$, so is $\Mod W \cap
Z$. It follows that $\Mod W \cap Z$ is an abelian category and that
the inclusion of it in $\Mod X$ is an exact functor commuting with
direct limits. It therefore has a right adjoint.
This proves (2).

Hence there is a commutative diagram  
\begin{equation}
\label{eq.W.cap.Z}
\begin{CD}
W \cap Z @>{\c}>> Z
\\
@V{\d}VV @VV{\b}V
\\
W @>>{\a}> X
\end{CD}
\end{equation}
of spaces and weak maps. 
Since $\a_*,\b_*,\c_*$ and $\d_*$, are inclusions of
subcategories, $\a_*\d_*=\b_*\c_*$. Up to natural equivalence we
also have $\d^!\a^!=\c^!\b^!$.

(3)
To see that $W \cap Z$ is closed in $W$ it suffices to show that 
$\d_*$ commutes with products; that is, that $\Mod W \cap Z$ is closed 
under products in
$\Mod W$. Let $M_\l$ be a family in $\Mod Z \cap W$, and suppose
they have a product in $\Mod W$. That product is $\a^!(\prod
M_\l)$, where $\prod M_\l$ is their product in $\Mod X$. But $\Mod
Z$ is closed under products in $\Mod X$, so $\prod M_\l$ belongs to
$\Mod Z$. Therefore its submodule $\a^!(\prod M_\l)$ is a 
$Z$-module, and hence
is in $\Mod W \cap Z$. This proves that $\d_*$ has a left adjoint
$\d^*$, whence $W \cap Z$ is closed in $W$.

We will now show that 
\begin{equation}
\label{eq.stars}
\c_*\d^* \cong \b^*\a_*.
\end{equation}
Since $\c_*\d^*$ is a left adjoint to $\d_*\c^!$ and $\b^*\a_*$ is
a left adjoint to $\a^!\b_*$, it suffices to show that $\d_*\c^!
\cong \a^!\b_*$. If $M$ is a $Z$-module, $\d_*\c^!M$ is the largest
submodule of $M$ that is in $\Mod Z \cap \Mod W$; since every
submodule of $M$ is in $\Mod Z$, $\d_*\c^!M$ is the largest
submodule of $M$ that is in $\Mod W$; but that submodule is
$\a^!\b_*M$. It follows that 
\begin{equation}
\label{eq.useful}
\d_*\c^! \cong \a^!\b_*,
\end{equation}
and (\ref{eq.stars}) follows from this.

(4)
If $W$ and $Z$ are closed subspaces of $X$, then $\Mod W$ and
$\Mod Z$ are closed under products in $\Mod X$. It follows that
$\Mod W \cap Z$ is closed under products in $X$, whence $W \cap Z$
is closed.

(5) and (6) follow from (2) and (3) because $Z \subset W$ implies
that $W \cap Z=Z$.

(7)
Since $Z$ is both weakly closed and weakly open in $X$, $\b_*$ has
an exact left adjoint $\b^*$. In particular, $Z$ is closed in $X$.
Hence by (3), $W \cap Z$ is closed in $W$, and $\c_*\d^* \cong
\b^*\a_*$. The right-hand side of this is a composition of exact
functors, so $\c_*\d^*$ is exact. It follows that $\d^*$ is exact
because $\c$ is a weakly closed immersion.
\end{pf}

{\bf Remarks. 1.}
The following two observations are easily checked.
If $W$ is a weakly closed subspace of $Z$, and $Z$ is a weakly
closed subspace of $X$, then $W$ is a weakly closed subspace of
$X$.  If $W$ is a closed subspace of $Z$, and $Z$ is a 
closed subspace of $X$, then $W$ is a closed subspace of
$X$.

{\bf 2.}
The binary operation $\cap$ on weakly closed subspaces is
idempotent, commutative, and associative.
One has $W \cap \phi=\phi$ and $W \cap X=W$.

{\bf 3.}
One may define the intersection of an arbitrary collection of
weakly closed subspaces $W_i$ by $\Mod (\cap W_i):=\cap \Mod W_i$.
It is clear that this is closed in $\Mod X$ 
under subquotients and direct limits because each $\Mod W_i$ is.
Hence $\cap_i W_i$ is a weakly closed subspace of $X$ (cf.
\cite[Lemma 6.2.2, Chapter 3]{Rosen}).

\begin{definition}
\label{defn.union}
Let $W_i$, $i \in I$, be a finite collection of weakly 
closed subspaces of
$X$. Their {\sf union}, denoted $\cup_i W_i$, is the smallest 
weakly closed closed subspace of $X$ that contains all the $W_i$.
It is the intersection of all the weakly closed subspaces that
contain all the $W_i$.
\end{definition}

{\bf Remarks. 1.}
We could first have defined the union of a pair of weakly closed
subspaces, thus making $\cup$ a binary operation. One sees 
that $\cup$ is idempotent, commutative, and associative. 
Associativity could then be used to define arbitrary finite unions;
this definition would agree with the one in Definition
\ref{defn.union}. It is clear that $W \cup \phi=W$ and $W \cup
X=X$. 

{\bf 2.}
If $W$ and $Z$ are weakly closed, then $\Mod W \cup Z$ is the 
full subcategory consisting of
all subquotients of modules of the form $M \oplus N$ where $M \in
\Mod W$ and $N \in \Mod Z$. Obviously this subcategory
is closed under subquotients; it is closed under direct sums
because $\Mod W$ and $\Mod Z$ are, and because 
direct sums are exact in $\Mod X$; it is therefore closed under
direct limits, so is weakly closed; it contains $\Mod W$ and
$\Mod Z$, and any weakly closed subcategory containing $\Mod W$ and
$\Mod Z$ must contain all modules of the form $M \oplus N$; hence this
category equals  $\Mod W \cup Z$.

{\bf 3.}
If $X$ is the affine space with coordinate ring $R$, and $W$ and $Z$
are the closed subspaces cut out the two-sided ideals $I$ and $J$,
then $W \cup Z$ is the closed subspace cut out by $I \cap J$. To see
this first observe that $\Mod R/I \cap J$ is weakly closed and
contains both $\Mod R/I$ and $\Mod R/J$; secondly, any weakly closed 
subspace that contains both $W$ and $Z$ must contain $R/I \oplus
R/J$ so must contain its submodule $R/I \cap J$, and therefore must
contain $\Mod R/I \cap J$.

{\bf 4.}
Let $W$ and $Z$ be weakly closed subspaces of $X$. Their {\sf
Gabriel product} $W \bullet Z$ is the weakly closed subspace
defined by declaring $\Mod W \bullet Z$ to be the full subcategory
of $\Mod X$ consisting of the modules $M$ for which there is an
exact sequence $0 \to L \to M \to N \to 0$ with $N \in \Mod W$ and
$L \in \Mod Z$ \cite[p. 395]{G}. Van den Bergh has shown that $W
\bullet Z$ is weakly closed, and is closed if both $Z$ and $W$ are
\cite[Prop. 3.4.5]{vdB}. Since $W \bullet Z$ contains both $W$ and
$Z$, it follows that 
\begin{equation}
\label{eq.cup}
W \cup Z \subset W \bullet Z \cap Z \bullet W.
\end{equation}
In the situation of the previous remark, $\Mod(W \bullet Z \cap Z
\bullet W) = \Mod R/IJ+JI$, so in general the inclusion in
(\ref{eq.cup}) is strict. Example \ref{eg.bad.triang} exhibits
closed points $p$ and $q$ in a an affine space such that $p \cup q=
p\bullet q \ne q \bullet p$.

{\bf 5.}
If $W,Y$ and $Z$ are weakly closed 
subspaces of $X$, we do not know whether
$
W \cap (Y \cup Z)= (W\cap Y) \cup (W \cap Z).
$
The analogous equality for weakly open subspaces does hold (Lemma 
\ref{lem.dist}). 

\section{Closed subspaces}
\label{sect.closed}

The first evidence that the notion of a closed subspace in
Definition \ref{defn.weakly.closed} is appropriate for
non-commutative geometry is the following result of Rosenberg
\cite[Proposition 6.4.1, p.127]{Rosen}. 
If $X$ is an affine space, say $X=\Mod R$, then each two-sided
ideal $I$ cuts out a closed subspace, namely $\Mod R/I$, and this
sets up a bijection between the closed subspaces of $X$ and the 
two-sided ideals of $R$ (see also \cite[Proposition 2.3]{SZ}). 
We call the closed subspace $\Mod R/I$ the {\sf zero locus} of $I$.
Thus, if $X$ is an affine scheme, its closed subschemes are in
natural bijection with the closed subspaces.
We prefer to say that the closed subspaces are the same things as
the closed subschemes.
Further evidence is provided by \cite[Theorem 3.2]{SmMaps} which
shows that a two-sided ideal in an $\NN$-graded $k$-algebra $A$
cuts out a closed subspace of the projective space with homogeneous
coordinate ring $A$.
There is further evidence.
If $X$ is a scheme such that $\Qcoh X$ is a Grothendieck category
(for example, if $X$ is noetherian) and $Z$ is a closed
subscheme of $X$, then $\Qcoh Z$ is a closed subspace of $\Qcoh X$ 
in the sense of Definition \ref{defn.weakly.closed}. Indeed, this
is the motivation for the definition of a closed subspace.

The main result in this section (Theorem \ref{thm.closed=closed}
shows that if  $X$ is a  quasi-projective scheme over a field $k$,
then the closed subschemes give all the closed subspaces.

We also prove a non-commutative geometric version of Matlis' result 
on the bijection between prime ideals in a commutative noetherian ring
and indecomposable injectives over that ring (Proposition
\ref{prop.bijection}).

\medskip

{\bf Remarks.}
{\bf 1.}
An interesting and ubiquitous non-commutative space is the graded
line $\LL^1$, defined by $\Mod \LL^1=\GrMod k[x]$ where $k$ is a
field and $\deg x=1$.  We associate to every subset 
$D \subset \ZZ$ a closed subspace $Z_D$ of $\LL^1$ defined by
declaring $\Mod Z_D$ to be the full subcategory consisting of those
modules $M$ such that $M_n=0$ for all $n \notin D$. 
Since there is a proper descending chain of 
subsets $D$, $\LL^1$ does not have the descending chain condition
on closed subspaces. If $D$ consists
of $n$ consecutive integers, then $Z_D$ is isomorphic to the affine
space with coordinate ring the $n \times n$ upper triangular 
matrices over $k$. 

If $X$ is an affine space with coordinate ring the enveloping algebra 
of a non-abelian solvable Lie algebra over a field $k$ of
characteristic zero, then $X$ contains many
weakly closed subspaces isomorphic to $\LL^1$ \cite{SZ}. If $\tilde
X$ is the blowup of a non-commutative surface $X$ at a closed point
$p$ lying on an effective divisor $Y$ such that $\{\cO_p(nY) \; |
\; n \in \ZZ\}$ is infinite, then the exceptional curve on
$\tilde X$ is a weakly closed subspace of $\tilde X$ isomorphic to 
$\LL^1$ \cite{vdB}.

{\bf 2.} 
The category of right comodules over a coalgebra is a
Grothendieck category, so is a space. Let $C$ and $D$ be coalgebras
and let $X$ and $Y$ be the spaces defined by $\Mod X=\Comod C$ and
$\Mod Y=\Comod D$. If $\varphi:C \to D$ is a
homomorphism of coalgebras, then the rule $(M,\rho) \mapsto
(M,(\id_M \otimes \varphi)\rho)$ gives an exact functor $f_*:\Comod
C \to \Comod D$. This is the direct image functor for an affine map
$f:X \to Y$. The right adjoint to $f_*$ is $f^!=- \square_D C$, the
co-tensor product. Joyal and Street \cite[Proposition 2, p.
454]{JS} show that the closed subspaces of $Y$ are in 
bijection with the subcoalgebras of $D$. Actually they don't quite
prove that, but using the fact that in a locally finite space
weakly closed and closed subspaces coincide (Proposition
\ref{prop.loc.finite}), their result implies this.
The closed points of $Y$, in the sense of Section 
\ref{sect.pts}, are in bijection with the simple subcoalgebras of $D$.

\begin{theorem}
\label{thm.closed=closed}
Let $X$ be a noetherian scheme having an ample line bundle (for
example, let $X$ be quasi-projective over $\Spec k$ with $k$ a
field). Then closed subschemes and closed subspaces of $X$ are the
same things.
\end{theorem}

We will prove the theorem after the next lemma.

Let $\cO_X(1)$ denote an ample line bundle on $X$.
If $i:Z \to X$ is the inclusion of a closed subscheme, then
the restriction $i^*(\cO_X(1))$ to $Z$ is ample.
We will denote it by $\cO_Z(1)$.

Let $\cM$ be an $\cO_X$-module. 
Its {\sf annihilator} is the kernel of the natural map 
$\cO_X \to \cHom_{\cO_X}(\cM,\cM)$. This is also the largest ideal
$\cI$ in $\cO_X$ such that $\cM\cI=0$, where $\cM\cI$ is
defined as the image of the natural map $\cM \otimes_{\cO_X} \cI
\to \cM$. We denote it by $\Ann \cM$.
If $\cM$ is coherent, its {\sf scheme-theoretic support} 
is the closed subscheme of $X$ cut out by its annihilator, i.e., 
it is $(Z,\cO_Z)$ where $\cO_Z=\cO_X/\Ann \cM$.

\begin{lemma}
\label{lem.inv.sub}
\label{lem.wkly.closed}
Let $\cM$ be a coherent $\cO_X$-module, and $Z$ its
scheme-theoretic support.
If $W$ is a weakly closed subspace of $X$ such that $\cM \in \Mod
W$, then $Z \subset W$, i.e., $\Qcoh Z \subset \Mod W$.
\end{lemma}
\begin{pf}
We show first that some
finite direct sum of copies of $\cM$ has a submodule that is
an invertible $\cO_Z$-module.

Since $\cM$ is a coherent $\cO_Z$-module, there is an epimorphism
$\bigoplus \cO_Z({-n}) \to \cM$ for some $n \gg 0$ and some suitable 
finite direct sum. 
The image of each $\cO_Z({-n})$ in $\cM$ is isomorphic
to $(\cO_Z/\cI_j)({-n})$ for some coherent ideal 
$\cI_j$ in $\cO_Z$.  Since $\cM$ is the sum of these images, 
the annihilator of $\cM$ in $\cO_Z$
is the intersection of all the $\cI_j$. By definition of $Z$,
that annihilator is zero,
so the diagonal map $\cO_Z \to \bigoplus \cO_Z/\cI_j$ is monic.
It follows that $\cO_Z({-n})$ embeds in a finite direct sum of copies 
of $\cM$.

By the preceding paragraph, $\cO_Z(-n)$ belongs to $\Mod W$ for 
some $n$, 
and hence for all $n \gg 0$. Since every coherent $\cO_Z$-module
is a quotient of a suitable finite direct sum of copies of $\cO_Z(-n)$
for some sufficiently large $n$, $\coh Z \subset \Mod W$. Since $\Mod W$
is closed under direct limits it follows that $\Qcoh Z \subset \Mod W$.
\end{pf}

\begin{pf} (of Theorem \ref{thm.closed=closed}.)
Let $W$ be a closed subspace of $X$. We must show that $\Mod W
=\Qcoh Z$ for some closed subscheme $Z$ of $X$. Let $\a:W \to
X$ be the inclusion. The counit $\cO_X
\to \a_*\a^*\cO_X$ is an epimorphism, so $\a_*\a^*\cO_X=\cO_Z$ for a 
unique closed subscheme $Z$ of $X$. Since $\cO_Z$ is a $W$-module,
Lemma \ref{lem.wkly.closed} shows that $\Qcoh Z \subset \Mod W$.

Suppose that $\cF \in \mod W$, and let $Y$ denote its
scheme-theoretic support.
By Lemma \ref{lem.wkly.closed}, $\Qcoh Y \subset \Mod W$. 
Since $Y$ is also the scheme theoretic support of
$\cF(n)$, $\cF(n) \in \Mod W$ for all $n \in \ZZ$.
For some $n \gg 0$ and some suitably large finite direct sum,
there is an epimorphism $\bigoplus \cO_X \to \cF(n)$, hence an
epimorphism
$$
\bigoplus \a_*\a^*\cO_X \to \a_*\a^*(\cF(n)) \cong \cF(n).
$$
It follows that $\cF(n)$ is in $\Qcoh Z$. Hence $\cF \in \Qcoh Z$
also, and we conclude that $\Mod W=\Qcoh Z$ as desired.
\end{pf}

This result suggests the following definitions.

\begin{definition}
Let $X$ be a locally noetherian space. The {\sf support} of a
noetherian $X$-module $M$ is the smallest closed subspace $Z$ such
that $M \in \Mod Z$. 
Such a smallest closed subspace exists because an intersection of
closed subspaces is closed.
A non-zero $X$-module is {\sf prime} if all
its non-zero submodules have the same support.
A non-empty closed subspace $Z$ of $X$ is {\sf prime} if it is the 
support of a prime noetherian $X$-module.
\end{definition}

The definition of a prime module emerged in discussions with C.
Pappacena and agrees with his definition in \cite{Papp}.

If $X$ has the descending chain condition on closed subspaces, 
then every non-zero $X$-module has a prime submodule: simply take
a submodule whose support is as small as possible. 
As remarked at the beginning of this section, the graded line,
$\GrMod k[x]$, does not have the descending chain condition
on closed subspaces; it is easy to see
that $k[x]$ does not contain a prime submodule.

If $X$ is affine with coordinate ring $R$, then the support of an
$R$-module is the zero locus of its annihilator because $R/\Ann M$
embeds in a product of copies of $M$, so is a module over any
closed subspace that contains $M$.
A closed subspace of $X$ is prime if and only if it is the zero
locus of a prime ideal; thus prime subspaces of $X$ 
are in bijection with the prime ideals in $R$.

A non-zero submodule of a prime module is prime.
A direct sum of copies of a prime module is prime.
An auto-equivalence of $\Mod X$ sends prime modules to prime
modules.

By Lemma \ref{lem.wkly.closed}, if $X$ is a noetherian scheme
having an ample line bundle, the support as we have just defined it
is the same as the scheme-theoretic support. 

\begin{proposition}
Let $X$ be a noetherian scheme having an ample line bundle.
\begin{enumerate}
\item{}
The (scheme-theoretic) support of a prime module is  integral.
 \item{}
A closed subscheme $Z \subset X$ is integral if and only if $\cO_Z$
is prime.
\end{enumerate}
\end{proposition}
\begin{pf}
We will use the fact that a closed subscheme
$Z$ of $X$ is reduced and irreducible if and only if the
product of non-zero ideals of $\cO_Z$ is non-zero.

(1)
Let $Z$ denote the support (= scheme-theoretic support) of a prime
module $\cM$. Without loss of generality we may assume that $\cM$
is coherent. By the arguments in Lemma \ref{lem.wkly.closed},
$\cO_Z(-n)$ embeds in a finite direct sum of copies of $\cM$. 
It follows that $\cO_Z(-n)$, and hence $\cO_Z$, is prime. 
Hence, the scheme-theoretic support  of every  non-zero coherent
ideal of $\cO_Z$ is $Z$. In other words, if $\cI$ is a non-zero 
coherent ideal of $\cO_Z$, then $\cI\cJ\ne 0$ for all
non-zero ideals $\cJ$ of $\cO_Z$. Hence $Z$ is reduced and
irreducible.

(2)
If $Z$ is integral, every non-zero ideal of $\cO_Z$ has annihilator
equal to zero, so has scheme-theoretic support $Z$, 
thus showing that $\cO_Z$ is prime.
Conversely, if $\cO_Z$ is prime, then every non-zero ideal of
it has support equal to $Z$ so has zero annihilator, thus showing
that $Z$ is integral.
\end{pf}

The conclusion of Lemma \ref{lem.wkly.closed} 
captures one of the essential aspects of commutativity:
it says that $X$ has a lot of closed subspaces.
This is similar to saying that a ring has lots of two-sided ideals. 
The next definition singles out this property and we think of it as
saying that the space is ``almost commutative''.

\begin{definition}
\label{defn.enough.closed}
We say that a locally noetherian space $X$ has {\sf enough closed
subspaces} if every noetherian module $M$ has the following
property: the smallest weakly closed subspace $Z \subset X$ such
that $M \in \Mod Z$ is equal to the support of $M$.
\end{definition}

In other words, $X$ has enough closed subspaces if every weakly 
closed subspace $W \subset X$ has the following property: 
if $M$ is a noetherian $W$-module, then there is a closed subspace
$Z$ of $X$ such that $Z \subset W$ and $M \in \Mod Z$.
Lemma \ref{lem.wkly.closed} implies that a quasi-projective scheme
has enough closed subspaces.

The next result shows that this notion of
``almost commutative'' is compatible with an older idea of what it
means for a ring to be ``almost commutative''. 
Gabriel's condition (H) \cite[page 422]{G} is expressed as part (2)
of the next result; Gabriel (loc. cit., page 423) 
showed that (2) implies (3) and Cauchon \cite{C} proved the
converse. The equivalent condition (1) is geometric
and therefore has meaning for non-affine spaces.
For example, if $A$ is a connected graded noetherian algebra that
is a finite module over its center, then property (1) is satisfied
by the projective space $X$ having $A$ as a homogeneous coordinate
ring.

\begin{proposition}
\label{prop.fbn}
Let $X$ be a locally noetherian affine space with coordinate ring
$R$. The following are equivalent:
\begin{enumerate}
\item{}
$X$ has enough closed subspaces;
\item{}
if $M$ is a noetherian $R$-module, then $\Ann M = \bigcap_{i=1}^n
\Ann(m_i)$ for some finite set of elements $m_i \in M$;
intersection
\item{}
there is a bijection between the prime ideals in $R$ and the
indecomposable injective $R$-modules.
\end{enumerate}
\end{proposition}
\begin{pf}
Gabriel and Cauchon have proved the equivalence of (2) and (3).

(2) $\Rightarrow$ (1)
Let $M$ be a noetherian $X$-module, and $W$ a weakly closed
subspace containing $M$.
Then $W$ contains $Z=\Mod R/\Ann M$.
There eare elements $m_i \in M$ such that
the diagonal map $R/\Ann M \to \bigoplus_{i=1}^n m_iR$ is injective.
Since $R/\Ann M$ embeds in a finite direct sum of copies of $M$,
whence $R/\Ann M \in \Mod W$, whence $Z \subset W$.

(1) $\Rightarrow$ (2)
Let $M$ be a noetherian $R$-module with support $Z$.
The smallest weakly closed subspace of $X$ containing $M$ is $W$,
where $\Mod W$ consists of all subquotients of direct sums
of copies of $M$. 
Then $R/\Ann M$ is in $\Mod W$, so is a subquotient of a
direct sum of copies of $M$. Since $R/\Ann M$ is a projective
module over itself it follows that  $R/\Ann M$ is a submodule of a
finite direct sum of copies of $M$. If $\varphi:R/\Ann M \to
M^{\oplus n}$ is injective, and $m_i$ denotes the image of $\varphi(1)$
under the $i^{\th}$ projection, then $\Ann M= \bigcap_{i=1}^n \Ann
(m_i)$. Thus (2) holds.
\end{pf}

Rings satisfying (2) are said to be fully bounded---see
\cite[Sections 6.4 and 6.10.4]{MR}. The standard example of such a
ring is one that is a finite module over its center.

In the next proof we use the notion of Krull dimension as defined
by Gabriel. A noetherian $X$-module is {\sf critical} if $\Kdim M/N
< \Kdim M$ for all non-zero submodules $N \subset M$. The
definition of Krull dimension is such that every non-zero 
module has a non-zero critical submodule. A submodule of a critical
module is critical.
Hence a prime subspace is the support of a critical 
module, and that module may also be chosen to be prime and
noetherian. 

\begin{proposition}
\label{prop.bijection}
Let $X$ be a locally noetherian space having the descending chain
condition on closed subspaces. 
If $X$ has enough closed subspaces there is
a bijection 
$$
\{\hbox{prime subspaces of $X$}\}
\longleftrightarrow
\{\hbox{indecomposable injective $X$-modules}\}.
$$
\end{proposition}
\begin{pf}
Let $Z$ be a prime subspace of $X$ and let $M$ be a critical
noetherian prime module having support equal to $Z$. Since $M$ is
critical it is uniform, so its injective envelope, say $E$, is 
indecomposable. We will show that the rule $Z \mapsto E$ is a
well-defined map setting up the claimed bijection.

To show this is well-defined, we must show that $E \cong E'$
whenever $E'$ is an injective envelope of another
critical noetherian prime module, say $M'$, whose support is
equal to $Z$. 

The hypothesis on $X$ implies that $\Mod Z$
consists of all subquotients of all direct sums of copies of $M$.
Of course, the same is true with $M'$ in place of $M$. Since each of
$M$ and $M'$ is a subquotient of a finite direct sum of copies of the
other they have the same Krull dimension. 

Choose the smallest integer $n$
such that $M'$ is isomorphic to $L/K$ where 
$K \subset L \subset M_1 \oplus \ldots \oplus M_n$ and each $M_i
\cong M$. The minimality of $n$ ensures that $L \cap M_i$ is
non-zero for all $i$.
If $K \cap M_i$ were non-zero for all $i$, then the Krull dimension
of $L/K$ would be strictly smaller that that of $M$; this is
not the case so $K \cap M_i$ is zero for some $i$. Hence $L \cap M_i$
embeds in $M'$. In particular, $M$ and $M'$ have a common non-zero
submodule, from which we conclude that $E \cong E'$.
Hence the map $Z \mapsto E$ is well-defined.

To see that this map is injective, suppose that $Z$ and $Z'$ are prime
subspaces for which the corresponding indecomposable injectives,
say $E$ and $E'$, are isomorphic. Let $M$ and $M'$ be  critical
noetherian prime modules whose supports are $Z$ and $Z'$
respectively. Thus $E$ and $E'$ are injective envelopes of $M$ and
$M'$ respectively. Since $E$ and $E'$ are isomorphic and
indecomposable, $M$ and $M'$ have a common non-zero submodule, say
$N$. However, since $M$ and $M'$ are prime they have the same
support as $N$. Thus $Z=Z'$.

To show that the map is surjective, suppose that $E$ is an
indecomposable injective. Because $X$ has the descending chain 
condition on closed subspaces, $E$ has a prime submodule and hence a 
prime critical noetherian submodule, say $M$. If $Z$ is the support
of $M$, then the map defined above sends $Z$ to $E$.
\end{pf}

The next result gives another class of spaces which have enough
closed subspaces. It applies for example to the category of
comodules over a coalgebra that is defined over a field.

\begin{proposition}
\label{prop.loc.finite}
Every weakly closed subspace of a locally finite space is closed.
\end{proposition}
\begin{pf}
Let $X$ be a locally finite space; that is, every $X$-module is a
direct limit of finite length $X$-modules. 
Let $W$ be a weakly closed subspace of $X$.
It suffices to show that $\Mod W$ is closed under products because
then the inclusion $i_*:\Mod W \to \Mod X$ will have a left
adjoint.

Let $\{M_i \; | \; i \in I\}$ be a collection of $W$-modules. By
hypothesis, $\prod_{i \in I} M_i$, their product in $\Mod X$, is a
sum of finite length $X$-modules. Let $N$ be a finite length
submodule of $\prod_{i \in I} M_i$. For each $j \in I$ let
$\rho_j:N \to M_j$ be the composition of the inclusion $N \to
\prod_{i \in I} M_i$ with the projection $\prod_{i \in I} M_i \to
M_j$. Write $K_j=\ker \rho_j$. Since $\cap_{j \in I} K_j$ is the
kernel of the inclusion $N \to \prod_{i \in I} M_i$ it is zero;
since $N$ has finite length there is a finite subset $J \subset I$
such that $\cap_{j \in J} K_j=0$. It follows that the composition
$$
N \to  \prod_{i \in I} M_i \to \prod_{j \in J} M_j
$$
is monic. Since $J$ is finite, $\prod_{j \in J} M_j$ is isomorphic
to $\oplus_{j \in J} M_j$, and therefore belongs to $\Mod W$. 
Hence $N$ belongs to $\Mod W$. 

We have shown that every finite length submodule of  $\prod_{i \in
I} M_i$ belongs to $\Mod W$. Since $\Mod W$ is closed under direct
limits we conclude that  $\prod_{i \in I} M_i$ belongs to $\Mod W$.
\end{pf}

\section{Closed points}
\label{sect.pts}

{\bf Notation.}
If $D$ is a division ring, we write $\Spec D$ for the space $\Mod D$ of
right $D$-modules.

\begin{definition}
\label{defn.ratl.pts}
A {\sf closed point} of a space $X$ is a closed subspace that is
isomorphic to $\Spec D$ for some division ring $D$.
We call it a {\sf $D$-rational point} of $X$.
\end{definition}

We denote a closed point by a single letter, say $p$, and often
write $p \in X$ to indicate that $p$ is a closed point of $X$. 
Since $\Mod p$ is equivalent to the category of modules over a
division ring it contains a unique simple module up to isomorphism, 
and because $\Mod p$ is closed under subquotients that simple 
is simple as an $X$-module.  We denote any module in this
isomorphism class by $\cO_p$. Thus
$D$ is isomorphic to $\End_X\cO_p$, and
every $p$-module is a direct sum of copies of $\cO_p$.

Although a closed point determines a simple module, a simple module
need not determine a closed point.

We will say that an $X$-module $M$ is {\sf compact} if
$\Hom_X(M,-)$ commutes with direct sums. 
If $X$ is locally noetherian, this is equivalent to the 
condition that $\Hom(M,-)$ 
commute with direct limits, and to the condition that $M$ be
noetherian \cite[Section 5.8]{Pop}.

\begin{definition}
An $X$-module $S$ is {\sf tiny} if $\Hom_X(M,S)$ is a finitely
generated module over $\End_XS$ for all compact
$X$-modules $M$.
\end{definition}

\begin{lemma}
\label{lem.prod.tiny}
Let $S$ be a tiny simple $X$-module. Then every compact
submodule of a direct product of copies of $S$
is isomorphic to a finite direct sum of copies of $S$.
\end{lemma}
\begin{pf}
Let $D=\End_X S$.
Let $P$ be a direct product of copies of $S$, and $M$ a compact
submodule of $P$. Define $M_0=M$ and, if $M_{i-1}$ is non-zero, we
define $M_i$ inductively as follows. Since $M_{i-1} \ne
0$, there is projection $P \to S$ which does not vanish on
$M_{i-1}$. Let $\rho_i:M \to S$ be the restriction of that projection 
and define $M_i=M_{i-1} \cap \ker\rho_i$.
Since $\dim_D \Hom_X(M,S)< \infty$, there is a smallest integer $n$
such that $\rho_{n+1} = \sum_{i=1}^n \mu_i \rho_i$ for some 
$\mu_i \in D$. Therefore $\ker \rho_{n+1}$ contains 
$\ker\rho_1 \cap \cdots \cap \ker\rho_n$. It follows that 
$0=M_n= \ker\rho_1 \cap \cdots \cap \ker\rho_n$. Hence the map
$\rho_1 \oplus \ldots \oplus \rho_n:M \to S^{\oplus n}$ 
is an isomorphism.
\end{pf}

{\bf Notation.}
If $S$ is a simple $X$-module, we write $\Sum S$
for the full subcategory consisting of all modules that are
isomorphic to a direct sum of copies of $S$.
If $D$ denotes the endomorphism ring of $S$, then $\Sum S$ is
equivalent to $\Mod D$.

\smallskip
In a Grothendieck category a sum of simple modules is isomorphic to
a direct sum of simple modules. This is not true for $(\Mod
\ZZ)^{\op}$.

\begin{lemma}
\label{lem.ModS}
Let $S$ be a simple $X$-module. The inclusion
$i_*:\Sum S \to \Mod X$ has a right adjoint $i^!$ given by
$$
i^!N := \hbox{the sum of all submodules of $N$ that are isomorphic 
to $S$}.
$$
\end{lemma}
\begin{pf}
If $f:M \to N$ is a map, then $f(i^!M) \subset i^!N$, so $i^!$ can
be defined on morphisms by sending a morphism to its restriction.
Thus $i^!$ really is a functor, and it takes values in $\Sum S$
because a sum of simple modules is a direct sum. It is clear that 
$\Hom_X(S,M)=\Hom_X(S,i^!M)$, so $i^!$ is right adjoint to $i_*$.

Alternatively, since $\Sum S$ is closed under subquotients, $i_*$
is exact, and since $\Sum S$ is closed under direct sums $i_*$
commutes with direct sums, whence $i_*$ has a right adjoint.
\end{pf}

\begin{theorem}
\label{thm.tiny.simple.point}
Suppose that $X$ is locally noetherian.
Let $S$ be a simple $X$-module. The following are equivalent:
\begin{enumerate}
\item{}
there is a closed point $p \in X$ such that $S \cong \cO_p$;
\item{}
$i_*:\Sum S \to \Mod X$ has a left adjoint;
\item{}
every direct product of copies of $S$ is isomorphic to a 
direct sum of copies of $S$;
\item{}
$S$ is tiny.
\end{enumerate}
\end{theorem}
\begin{pf}
(1) $\Rightarrow$ (2)
Since $\Mod p$ consists of all direct sums of copies of $\cO_p$,
$\Mod p=\Sum S$. Hence the fact that $p$ is closed ensures that
$i_*$ has a left adjoint.

(2) $\Rightarrow$ (1)
By Lemma \ref{lem.ModS}, $i_*$ has a right adjoint, so
the hypothesis ensures that $\Sum S$ satisfies the requirements
to be the module category of a closed point.

(2) $\Leftrightarrow$ (3)
Since $i_*$ is exact, it has a left adjoint if and only if it
commutes with products. That is, if and only if $\Sum S$ is closed
under products. But this is equivalent to condition (3).

(3) $\Rightarrow$ (4)
If $M$ is compact, then $\Hom_X(M,S)=\Hom_X(M,i_*S) \cong
\Hom_X(i^*M,S)$. But $i^*M$ is a quotient of $M$, so is also of
compact.
But it is also in $\Sum S$, so is a {\it finite} direct sum of 
copies of $S$. Thus $\Hom_X(M,S)$ is a finite direct sum of copies
of $\End_X S$.

(4) $\Rightarrow$ (3)
Let $P$ be a product of copies of $S$. Since $X$ is locally noetherian,
$P$ is a direct limit of noetherian submodules. Each of those 
submodules is of
compact, and hence semisimple by Lemma \ref{lem.prod.tiny}.
Therefore $P$ is semisimple.
\end{pf}

A locally noetherian space $X$ over a field $k$ is {\sf Ext-finite}
if $\dim_k \Ext_X^i(M,N) < \infty$ for all noetherian $X$-modules
$M$ and $N$ and all integers $i$.
If $A$ is a connected graded $k$-algebra such that $\dim_k A_n <
\infty$ for all $n$, then \cite[Section 3]{AZ} gives conditions
which ensure that the projective space with homogeneous coordinate
ring $A$ is Ext-finite. 

\begin{corollary}
\label{cor.one.pt}
Let  $X$ be a locally noetherian space over $\Spec k$. If $X$ is
Ext-finite, then every simple $X$-module is isomorphic to $\cO_p$ for
some closed point $p$. If $k$ is algebraically closed, then $p$ is
$k$-rational.
\end{corollary}
\begin{pf}
Let $S$ be a simple $X$-module. Then 
$\Hom_X(M,S)$ is finite dimensional over $k$ whenever 
$M$ is noetherian, so $S$ is tiny. Hence there is a closed
point $p \in X$ such that $\cO_p \cong S$. Since $X$ is Ext-finite,
$\Hom_X(\cO_p,\cO_p)$ is a finite dimensional division algebra over
$k$, and therefore equal to $k$ if $k$ is algebraically closed.
\end{pf}

The following is routine but we record it for later reference.

\begin{proposition}
\label{prop.affine.points}
Let $X$ be an affine space with coordinate ring $R$. There
is a bijection between
\begin{enumerate}
\item{}
closed points in $X$;
\item{}
simple $R$-modules $S$ such that $S$ has finite length over $\End_R S$;
\item{}
maximal two-sided ideals $\fm$ such that $R/\fm$ is artinian.
\end{enumerate}
If $R$ is a countably generated algebra over an
uncountable algebraically closed field $k$, then all closed points are 
$k$-rational, and they are in bijection with the finite dimensional 
simple $R$-modules.
\end{proposition}

The hypotheses on $k$ and $R$ in the last part of the proposition are needed
because it is an open problem whether there exists a division algebra $D$
over an algebraically closed field $k$ such that $D \ne k$ but $D$
is finitely generated as a $k$-algebra. If there were such a $D$, then 
$\Mod D$ would have a single closed point but that point would not be 
$k$-rational and would not correspond to a simple module of finite 
dimension over $k$.


The next result says that point modules over a reasonable 
graded algebra yield closed points in the associated projective space.

\begin{proposition}
\label{prop.proj.pts}
Let $A$ be a right noetherian locally finite $\NN$-graded
$k$-algebra.  Let $X$ be the projective space with homogeneous 
coordinate ring $A$. 
Let $M$ be an infinite dimensional graded right $A$-module such 
that $\dim_k M/N <\infty$
for all non-zero graded submodules $N \subset M$.
If $A$ satisfies the condition $\chi_1$ in \cite[Section 3]{AZ},
then $\pi M \cong \cO_p$ for some closed point $p \in X$.
\end{proposition}
\begin{pf}
The $\chi_1$ hypothesis implies that 
$\Hom_X(\cO_X, \pi M(n))$ has finite dimension for all $n \in \ZZ$;
The hypothesis on $M/N$ ensures that $\pi M$ is a simple $X$-module. 
If $\cN$ is any right noetherian $X$-module, then there is an
epimorphism from a finite direct sum of copies of $\cO_X(n)$ for
various integers $n$ to $\cN$. Hence $\Hom_X(\cN,\pi M)$ is a
subspace of a finite direct sum of various  $\Hom_X(\cO_X(n), \pi
M)$, so finite dimensional. Thus $\pi M$ is a tiny simple.
\end{pf}

We have the following consequences of Proposition \ref{prop.W.cap.Z}.
If $W$ is a weakly closed subspace of $X$, and $p$ is a closed 
point of $X$ such that $\cO_p$ is a $W$-module, then $p$ is 
a closed point of $W$. 
If $Z$ is a closed subspace of $X$ and $p$ is a closed point of $Z$, 
then $p$ is a closed point of $X$. 
However, a closed point of a weakly closed subspace of $X$
need not be a closed point of $X$; for example, if
$M$ is an infinite dimensional simple module over a ring $R$ such
that $R/\Ann M$ is not artinian, then $M$ does not correspond to a
closed point of $\Mod R$ even though $\Sum M$ is a closed point in
$\Sum M$.

\smallskip

{\bf Remark.}
Closed points do not always behave like zero-dimensional
geometric objects. The following example has been emphasized by
Artin and Zhang.

Fix an algebraically closed field $k$.
Let $E$ be an elliptic curve over $\Spec k$, $\s$ an
automorphism of $E$ having infinite order, and $\cL$ a line bundle
of degree $\ge 3$ on $E$. Let $B=B(E,\s,\cL)$ be the twisted 
homogenous coordinate ring associated to this data \cite{ATV1}.
Then $B$ is a connected graded $k$-algebra of Gelfand-Kirillov
dimension two, so $A:=B \otimes_k B$ is the homogeneous coordinate
ring of a non-commutative projective space $X$ that is like a
non-commutative 3-fold. If $\cF$ and $\cG$ are noetherian
$X$-modules then $\Ext_X^i(\cF,\cG)$ is finite dimensional for all
$i$. The ring $A$ is Artin-Schelter-Gorenstein so $X$
satisfies Serre duality.
We use the $B$-$B$-bimodule structure 
on $B$ to view it as a graded right $A$-module.
Applying the quotient functor $\pi:\GrMod A \to \Tails
A =\Mod X$ to this module produces an $X$-module $\cM$.
Since $\s$ has infinite order every non-zero two-sided graded 
ideal in $B$ has finite codimension \cite{ATV1}; hence
$\cM$ is a simple $X$-module. By Corollary \ref{cor.one.pt}, there
is a closed $k$-rational point $p \in X$ such that $\cM \cong
\cO_p$. Because the Gelfand-Kirillov dimension of $B$ is two, and
$A$ is Cohen-Macaulay with respect to GK-dimension,
$\uExt^2_{A}(B,A) \ne 0$. By \cite[Theorem 8.1(5)]{AZ},
$\uExt^2_X(\cO_p,\cO_X) \ne 0$. Serre duality now gives
$H^1(X,\cO_p(i)) \ne 0$ for some integer $i$.
The non-vanishing of this cohomology group indicates that $p$ is
behaving like a geometric object of dimension greater than zero.

Another example of this sort of phenomenon is given in \cite{JS1}.
Let $k$ be the complex field.
There is a smooth non-commutative projective surface $X$ 
satisfying Serre
duality, and a closed $k$-rational point $p \in X$ such that
$\Ext^1(\cO_p,\cO_p)=0$ and $\Ext^2_X(\cO_p,\cO_p)=k$.
Thus $p$ behaves like a $-2$-curve on $X$. This analogy is
strengthened by the fact that $X$ is constructed from the sheaf of
differential operators on $\PP^1$ (and the projectivized 
cotangent bundle to $\PP^1$ is a surface with a $-2$-curve). 
An alternative construction of $X$ is to take a 
projective compactification (as
described before Proposition \ref{prop.Rees.pts}) of the affine
space with coordinate ring a non-simple primitive factor ring of
the enveloping algebra
$U(\fsl_2)$; the associated graded ring of this factor is the
coordinate ring of a singular quadric surface in $\AA^3$, and 
the exceptional fiber in
its minimal resolution is a $-2$-curve.
All this suggests that $p$ behaves both like a point and 
like a $-2$-curve. Van den Bergh has suggested that we should
perhaps think of $X$ as its own blowup at $p$. The plausibility of
this idea is reinforced by the fact that the ideal annihilating the
simple $\fsl_2$-module corresponding to $p$ is idempotent (because
$\fsl_2$ is semisimple), so the Rees ring construction of the
blowup would simply reproduce the original surface.

\begin{definition}
\label{defn.pt.image}
Let $f:Z \to X$ be a weak map of spaces.
Let $W$ be a weakly closed or weakly open subspace of $Z$ and let
$p$ be a closed point of $X$. We say that the {\sf image of $W$ is
$p$}, and write $f(W)=p$, if $0 \ne f_*(\Mod W) \subset \Mod p$.
\end{definition}

In that case there is a weak map $g:W \to p$ such that $ig=f\a$,
where $\a:W \to Z$ and $i:p \to X$ are the inclusions; its direct
image functor is $g_*=i^*f_*\a_*$.

The map $f:\Spec k \to X$ induced by the inclusion
$$
\begin{pmatrix} k & k \\ 0 & k \end{pmatrix}
\to
\begin{pmatrix} k & k \\ k & k \end{pmatrix}
$$
does not send the closed point of $\Spec k$ to a closed point of
$X$. The problem is that $f_*(\cO_{\Spec k})$ is not a direct sum
of copies of a single simple $X$-module. Similarly, if $f:\Spec
k \to \Spec (k \times k) $ is the map induced by the inclusion of
the diagonal matrices in $M_2(k)$, $f$ does not send the closed
point of $\Spec k$ to a closed point of $\Spec (k \times k) $.

\medskip
{\bf The bimodule $o_p$.}
Let $i:p \to X$ be the inclusion of a $D$-rational closed point of $X$.
Following \cite[Chapter 3]{vdB}, $p$ determines an $X$-$X$-bimodule 
$o_p$ defined by
$$
\cHom_X(o_p,-) = i_*i^!= \hbox{the sum of all submodules that are
isomorphic to $\cO_p$}.
$$
Its left adjoint is $-\otimes_X o_p=i_*i^*$.
If $\a:\Spec D \to p$ denotes the isomorphism
we may assume that $\a^*=\a^!$.
By definition $\Hom_X(\cO_p,-)=\a^!i^!$.
It is natural to denote the left adjoint of this by $-\otimes_D \cO_p$. We
write
$$
(-)^*=\Hom_D(-,D):\Mod D^{\op} \to \Mod D.
$$

The next result extends \cite[Lemma 5.5.1]{vdB}.

\begin{lemma}
\label{lem.vdB.4.6.1}
Let $p$ be a closed $D$-rational point of $X$.
For all $j \ge 0$,
\begin{enumerate}
\item{}
$\cExt^j_X(o_p,-) \cong \Ext^j_X(\cO_p,-) \otimes_D \cO_p$, and
\item{}
$\cTor^X_j(-,o_p) \cong \Ext^j_X(-,\cO_p)^* \otimes_D \cO_p$.
\end{enumerate}
\end{lemma}
\begin{pf}
(1)
Using the exactness of some of the functors, we have
\begin{align*}
\Ext^j_X(\cO_p,-) \otimes_D \cO_p &= (i_*\a_*) \circ R^j(\a^!i^!)
\\
&\cong R^j(i_*\a_*\a^!i^!)
\\
&\cong R^j(i_*i^!)
\\
&=\cExt^j_X(o_p,-).
\end{align*}

(2)
We have
$$
M \otimes o_p=i_*i^*M=i_*\a_*\a^*i^*M=(\a^*i^*M) \otimes_D \cO_p,
$$
and
\begin{align*}
\a^*i^* M & \cong \Hom_D(\Hom_D(\a^*i^*M,D),D)
\\
&
\cong \Hom_D(\Hom_X(M,i_*\a_*D),D)
\\
&= \Hom_X(M,\cO_p)^*.
\end{align*}
Combining these two computations shows that
$$
M \otimes_X o_p \cong \Hom_X(M,\cO_p)^*\otimes_D \cO_p.
$$
Taking left derived functors of both sides, and using the fact that $(-)^*$
and $-\otimes_D \cO_p$ are exact, gives (2).
\end{pf}

If $p$ and $q$ are distinct closed points, then $o_p \otimes_X o_q=0$.

\section{The complement to a weakly closed subspace}

Throughout this section $X$ is a locally noetherian space.

Weakly closed and weakly open subspaces were defined in section two.
We begin this section by defining the complement to a weakly closed
subspace. Results in Gabriel's thesis show that the complement to a
weakly closed subspace is weakly open, and that every weakly 
open subspace arises as such a complement.

\begin{definition}
Let $X$ be a locally noetherian space, and let
$W$ be a weakly closed subspace.
The category of $X$-modules {\sf supported on $W$} is the full
subcategory consisting of those $X$-modules which are
the directed union of a family of submodules $M$ each of which 
has a finite fitration by submodules $0=M_0 \subset M_1 \subset 
\ldots \subset M_r=M$ such that each $M_i/M_{i-1}$ is a $W$-module.
We denote this category by
$$
\Mod_W X.
$$
We call the modules $M_i/M_{i-1}$ the {\sf slices} of the
filtration.
\end{definition}

\begin{proposition}
\cite[Prop. 8, p. 377]{G}
Let $W$ be a weakly closed subspace of a locally noetherian space $X$. 
Then $\Mod_W X$ is a localizing subcategory of $\Mod X$.
\end{proposition}

\begin{definition}
\label{defn.sat}
The {\sf saturation} of a weakly closed subspace $W$ of $X$
is the weakly closed subspace $W_{sat}$ defined by
$\Mod W_{sat}= \Mod_W X$. We say that $W$ is {\sf saturated}
if $W=W_{sat}$.
\end{definition}

{\bf Remarks. 1.}
It is obvious that $(W_{sat})_{sat}=W_{sat}$.

{\bf 2.}
The proof of Lemma \ref{lem.union.open}
shows that if $W$ and $Z$ are weakly closed subspaces of a locally
noetherian space $X$, then $Z_{sat} \cap W_{sat} = (Z \cap
W)_{sat}$. 
This implies that $(Z \cap W_{sat})_{sat}= (Z \cap W)_{sat}$.

\begin{definition}
\label{defn.complement}
If $W$ is a weakly closed subspace of a locally noetherian space $X$, 
its {\sf complement}, denoted $X \backslash W$, 
is the non-commutative space defined by
\begin{equation}
\label{eq.open}
\Mod X \backslash W:= \Mod X /\Mod_W X.
\end{equation}
We write $j:X\backslash W \to X$ for the map consisting of the
adjoint pair of functors $(j^*,j_*)$ where 
$j^*:\Mod X \to \Mod X \backslash W$ 
is the quotient functor, and $j_*$ is its right adjoint (which
exists because $\Mod_W X$ is a localizing subcategory).
\end{definition}

One sees at once that $X \backslash X=\phi$ and $X \backslash \phi=X$.

\smallskip

Because $j_*$ is fully faithful we have the following result.

\begin{proposition}
Let $W$ be a weakly closed subspace of a locally noetherian space
$X$. Then $j:X \backslash W \to X$ is a weakly open immersion.
\end{proposition}

We will identify $\Mod X\backslash W$ with the full subcategory 
of $\Mod X$ consisting of those $X$-modules $M$ for which
the canonical map $M \to j_*j^*M$ is an isomorphism.
In this way we can speak of $X \backslash W$ as a 
weakly open subspace of $X$.

Since $\Mod_WX$ is a localizing subcategory of $\Mod X$, its
inclusion has a right adjoint $\tau:\Mod X \to \Mod_W X$. This
sends an $X$-module to its largest submodule that is supported on
$W$. We sometimes call it the {\sf $W$-torsion} functor.
By \cite{G}, for every $X$-module $M$ there is an exact sequence
\begin{equation}
\label{eq.tors.seq}
0 \to \tau M \to M \to j_*j^*M \to R^1\tau M \to 0,
\end{equation}
and $\tau j_*j^*=0$.

\begin{proposition}
Every weakly open subspace of $X$ is of the form $X \backslash W$
for some weakly closed subspace $W$.
\end{proposition}

If $Z$ is a closed subspace of
$X$, we call $X \backslash Z$ an {\sf open} subspace of $X$.

\medskip

The following two results are used in section \ref{sect.contains}.

\begin{lemma}
\label{lem.gabriel}
Let $W$ be a weakly closed subspace of $X$ and let $U=X\backslash
W$. Let $\a:W \to X$ and $j:U \to X$ be the inclusions. Then
$\Ext^1_X(\a_*M,j_*N)=0$ for all $M \in \Mod W$ and all $N \in \Mod
U$.
\end{lemma}
\begin{pf}
This result is a restatement of condition (b) in
\cite[Lemme 1, p. 370]{G}; condition (c) of that
result holds because $j^*j_* \cong \id_U$ and 
$j^*$ is a left adjoint to $j_*$. 
\end{pf}

\begin{proposition}
\label{prop.gabriel}
Let $U$ be a weakly open subspace of a locally noetherian space $X$. Let
$j:U \to X$ denote the inclusion. Then $j_*$ commutes with direct
limits.
\end{proposition}
\begin{pf}
This is a restatement of \cite[Coroll. 1, p. 379]{G}.
\end{pf}

\begin{definition}
\label{defn.union.open}
Let $U_i$, $i \in I$, be a small set of weakly open subspaces of
$X$. Let $j_i:U_i \to X$ be the inclusion maps. We define their {\sf
union}, $\cup_{i \in I} U_i$, to be the weakly open subspace of 
$X$ whose module category is $\Mod X/\sT$ where $\sT$ is 
the localizing subcategory consisting of the $X$-modules $M$ such 
that $j_i^*M=0$ for all $i$.
We say that the $U_i$ provide a {\sf weakly open cover} of $X$
if their union is $X$; i.e., if $\sT=0$. This is equivalent to the
requirement that the diagonal $\Mod X \to \prod_i \Mod U_i$ is
faithful. If each $U_i$ is open, we say that they provide
an {\sf open cover} of $X$.
\end{definition}

We could first define the union of a pair of weakly open
subspaces thus making $\cup$ a binary operation on weakly open
subspaces. This operation is idempotent, 
commutative, and associative. We can
use the associativity to define the union of a finite collection of
weakly open subspaces in the obvious way, and this agrees with the
definition above. One has $U \cup \phi=U$ and $U \cup X=X$.

\begin{lemma}
\label{lem.union.open}
Let $W$ and $Z$ be weakly closed subspaces of a locally noetherian
space $X$. If $U=X \backslash W$
and $V=X \backslash Z$, then $U \cup V =X \backslash (W \cap Z)$.
\end{lemma}
\begin{pf}
By definition, $\Mod U \cup V$ is the quotient of $\Mod X$  by
the localizing subcategory $\Mod_W X \cap \Mod_Z X$, so it suffices
to show that this is equal to $\Mod_{W \cap Z} X$. It suffices to
prove that the full subcategories of noetherian modules are the
same. It is obvious that $\mod_{W \cap Z} X \subset \mod_W X \cap
\mod_Z X$. 

To prove the reverse inclusion, suppose that $M$ is in
$\mod_W X \cap \mod_Z X$.
There are finite filtrations $0=L_0 \subset L_1 \subset \ldots
\subset L_r = M$ and $0=N_0 \subset N_1 \subset \ldots \subset
N_s=M$ such that $L_i/L_{i-1} \in \mod W$ and $N_j/N_{j-1} \in \mod
Z$ for all $i$ and $j$. Taking a common refinement of these 
filtrations produces a finite filtration of $M$ all of whose slices
belong to $\mod W \cap Z$, thus showing that $M \in \mod_{W \cap Z}
X$.
\end{pf}

{\bf Remarks. 1.}
It is straightforward to show that if $\sS \subset \sT \subset \sA$
are Serre subcategories of a Grothendieck category $\sA$, then
$\sT/\sS$ is a Serre subcategory of $\sA/\sS$, and $\sA/\sT$
is equivalent to $(\sA/\sS)/(\sT/\sS)$ \cite[Ex. 6, p. 174]{Pop}.

This has the following geometric interpretation. 
If $W \subset Z \subset X$ are
weakly closed subspaces of $X$, then $X\backslash Z$ is a weakly
open subspace of $X \backslash W$ (Corollary \ref{cor.open.in.open}).
Example \ref{eg.bad.triang} shows that $X \backslash Z$ is {\it not}
the complement to $Z\backslash W$ in $X \backslash W$; indeed, the
last phrase does not always make sense because $Z \backslash W$
need not be a weakly closed subspace of $X \backslash W$. 

{\bf 2.}
It follows from Lemma \ref{lem.union.open} and the previous 
remark that $U$ and $V$ are weakly open subspaces of $U \cup V$.
It is easy to see that $U \cup V$ is the smallest weakly open
subspace of $X$ that contains both $U$ and $V$.

\begin{example}
\label{eg.bad.triang}
There is a space $X$ with closed subspaces $W \subset Z \subset X$
for which there is no weak map from $Z\backslash W$ to $X \backslash W$
making the following diagram commute:
$$
\begin{CD}
Z\backslash W @. X \backslash W
\\
@VVV @VVV \\
Z @>>> X.
\end{CD}
$$
Let $k$ be a field, and let $X=\Mod R$ be the affine space 
with coordinate ring 
$$
R=\begin{pmatrix}
k & k 
\\
0 & k
\end{pmatrix}.
$$ 
Let $p$ and $q$ be the closed points of $X$
corresponding to the simple right modules
$\cO_p$ and $\cO_q$, labelled so that $\cO_p$ is the projective simple.
Let $Z=\Mod R/N$ where $N$ is the ideal of strictly upper 
triangular matrices,
and let $W=q$. Clearly $\Mod Z \cong \Mod k^{\oplus 2}$.
 Let $j:X \backslash W \to X$ be the inclusion. 
There is a non-split exact sequence 
$0 \to \cO_p \to j_*j^*\cO_p \to \cO_q \to 0$.
As full subcategories of $\Mod X$, $\Mod Z$ consists of all direct 
sums of copies of $\cO_p$ and $\cO_q$, $\Mod Z \backslash W$ consists of
all direct sums of copies of $\cO_p$, and $\Mod X \backslash W$ consists of
all direct sums of copies of the indecomposable module $j_*j^*\cO_p$.
The truth of the claim is now apparent.
\end{example}

This example is typical of what can happen in non-commutative geometry.
Any space having a non-split extension of 
two non-isomorphic tiny simple modules has a closed subspace 
isomorphic to $X$ (for example, the affine space with coordinate
ring the
enveloping algebra of any non-abelian solvable Lie algebra).
The behavior is due to the fact that 
$\Ext^1_X(\cO_q,\cO_p)\ne 0$ even though $p$
and $q$ are distinct closed points. 
Although $Z$ consists of two closed points
and $W$ consists of only one of them, 
$\Mod Z \cap \Mod X \backslash W=\phi$.

\begin{lemma}
Let $j:U \to X$ be the inclusion of a weakly open subspace. Then
\begin{enumerate}
\item{}
every weakly open subspace of $U$ is a weakly open subspace of $X$;
\item{}
if $i:V \to X$ is the inclusion of another weakly open subspace and
$\Mod V$ is contained in $\Mod U$, then there is a
weakly open immersion $\b:V \to U$ such that $i=j\b$.
\end{enumerate}
\end{lemma}
\begin{pf}
(1)
Let $\b:V \to U$ be the inclusion of a weakly open subspace of $U$.
The functor $i_*:=j_*\b_*$ is fully faithful, so we can view $\Mod
V$ as a full subcategory of $\Mod X$. And $i_*$ has an exact left
adjoint $\b^*j^*$, so $V$ is weakly open in $X$.

(2)
Let $\b_*:\Mod V \to \Mod U$ be the inclusion. Clearly
$i_*=j_*\b_*$. Since $j^*$ and $i^*$ have left adjoints, $\Mod U$
and $\Mod V$ are closed under products in $\Mod X$. It follows that
$\Mod V$ is closed under products in $\Mod U$. 
Since $j^*j_* \cong \id_U$, $\b_* \cong j^*i_*$, and this is 
left exact. Therefore $\b_*$ has a left adjoint, say $\b^*$. 
Notice that $i^* \cong \b^*j^*$. To see
that $\b^*$ is exact, suppose that $0 \to L \to M \to N \to 0$ is
exact in $\Mod U$. By \cite[Coroll. 1, p. 368]{G}, there is an
exact sequence $0 \to L' \to M' \to N' \to 0$ in $\Mod X$ such that
the first sequence is obtained by applying $j^*$ to the second. It
follows that $\b^*$ applied to the first sequence is isomorphic to
$i^*$ applied to the second; since $i^*$ is exact, we deduce that
$\b^*$ is exact. Thus $\b$ is the inclusion of an open subspace of
$U$.
\end{pf}

\begin{lemma}
\label{lem.int.opens}
Let $W$ and $Z$ be weakly closed subspaces of a locally noetherian
space $X$. Then
\begin{equation}
\label{eq.loc.cats}
\Mod_{W \cup Z} X=\Mod_{W\bullet Z} X= \Mod_{Z \bullet W} X,
\end{equation}
and this is the smallest localizing subcategory of $\Mod X$ that
contains $\Mod W$ and $\Mod Z$.
\end{lemma}
\begin{pf}
Since $W \cup Z$ contains $W$ and $Z$, 
$\Mod_{W \cup Z} X$ contains $\Mod W$ and $\Mod Z$. Since it is
closed under extensions it also contains $\Mod W \bullet Z$ and
$\Mod Z \bullet W$, and hence contains $\Mod_{W\bullet Z} X$ and
$\Mod_{Z \bullet W} X$. However, since $W \cup Z$ is contained in
$W\bullet Z$ and in $Z\bullet W$, $\Mod_{W \cup Z} X$ is contained
in  $\Mod_{W\bullet Z} X$ and  $\Mod_{Z \bullet W} X$.
This gives (\ref{eq.loc.cats}). Finally, any localizing category
that contains $\Mod W$ and $\Mod Z$ must contain $\Mod W \bullet
Z$ and therefore $\Mod_{W\bullet Z} X$, so this is the smallest
localizing subcategory that contains  $\Mod W$ and $\Mod Z$.
\end{pf}

\begin{definition}
\label{defn.intersect.opens}
Let $U$ and $V$ be weakly open subspaces of a locally noetherian
space $X$. The {\sf intersection} of $U$ and $V$ is
$$
U \cap V:=X \backslash (W \cup Z),
$$
where $W$ and $Z$ are weakly closed subspaces such that
$U=X \backslash W$ and $V=X \backslash Z$.
\end{definition}

{\bf Remarks. 1.}
The definition of $U \cap V$ does not depend on the choice of $W$
and $Z$.  Although $U$ may be the complement of many different 
weakly closed subspaces $W$, $\Mod_W X$ depends only on $U$
because, if $j:U \to X$ is the inclusion, $\Mod_W X$ consists of 
the modules on which $j^*$ vanishes;
by Lemma \ref{lem.int.opens}, $\Mod_{W \cup Z}X$ is 
the smallest localizing category containing $\Mod_W X$ and $\Mod_Z X$.

{\bf 2.}
It follows from the remark before Example \ref{eg.bad.triang}
that $U \cap V$ is a weakly open subspace of both $U$ and $V$.
Proposition \ref{prop.int.opens} shows that it is the largest 
weakly open subspace of $X$ that is contained in both $U$ and $V$.

{\bf 3.}
By Lemma \ref{lem.int.opens}, if $W$ and $Z$ are weakly closed
subspaces of $X$, and $U=X \backslash W$ and $V=X \backslash Z$ are
their complements, then 
\begin{equation}
\label{eq.int.opens}
U \cap V = X \backslash (W \bullet Z) = X\backslash (Z \bullet W).
\end{equation}

{\bf 4.}
If $W$ and $Z$ are closed subspaces of $X$, so is $Z \bullet W$
\cite[Prop. 3.4.5]{vdB}. Hence, if $U$ and $V$ are open subspaces
of $X$, so is $U \cap V$.

\begin{proposition}
\label{prop.int.opens}
Let $U$ and $V$ be weakly open subspaces of a locally noetherian
space $X$. Then $\Mod U \cap V = \Mod U \cap \Mod V.$
\end{proposition}
\begin{pf}
Let $W$ and $Z$ be weakly closed subspaces such that $U=X
\backslash W$ and $V=X \backslash Z$.
Let $i:U \to X$, $j:V \to X$, and $\a:U \cap V \to X$ be the
inclusions. To make sense of the statement of the proposition, we
identify $\Mod U$ with the full subcategory of $\Mod X$ consisting
of those modules $M$ for which the natural map $M \to i_*i^*M$ is
an isomorphism, etc.

The quotient functor $\a^*:\Mod X \to \Mod U \cap V$
vanishes on $\Mod_W X$, so $\a^*$ factors through $i^*:\Mod X \to
\Mod U$. Hence $\Mod U \cap V \subset \Mod U \cap \Mod V.$ 

It remains to show that if $M$ is in $\Mod U$ and $\Mod V$, then it
is in $\Mod U \cap V$. Let $\tau:\Mod X \to \Mod_{W \cup Z} X \to
\Mod X$ be the torsion functor. We must show that the middle arrow
in the exact sequence
$$
0 \to \tau M \to M \to \a_*\a^*M \to R^1\tau M \to 0
$$
is an isomorphism. If $\tau M$ is non-zero, then $M$ has a non-zero
submodule belonging to $\Mod W \cup Z$. By (\ref{eq.cup}), $M$
therefore has a non-zero submodule belonging to either $\Mod W$ or
$\Mod Z$. But $M \cong i_*i^*M \cong j_*j^*M$, so it cannot have
such a non-zero submodule; we conclude that $\tau M=0$. If $R^1\tau
M \ne 0$, then there is an essential extension $0 \to M \to M' \to
T \to 0$ with $T \in \Mod W \cup Z$. By (\ref{eq.cup}), $T$ has a
non-zero submodule belonging to either $\Mod W$ or $\Mod Z$, say
$T'$, and this would give an essential extension $0 \to M \to M''
\to T' \to 0$ contradicting the fact that  $M \cong i_*i^*M \cong
j_*j^*M$. Thus $R^1\tau M=0$ also, and the result follows.
\end{pf}

{\bf Remarks. 1.}
It follows from Proposition \ref{prop.int.opens} that $U \cap V$ is
the largest weakly open subspace of $X$ contained in both $U$ and $V$.

{\bf 2.}
The binary operation $\cap$ on weakly open subspaces is idempotent,
commutative, and associative. The associativity allows us to define
the intersection of a finite number of weakly open subspaces by
induction.
Both $\phi$ and $X$ are open subspaces of $X$, and if $U$ is weakly
open, then $U \cap \phi=\phi$ and $U \cap X=X$.
We show next that the lattice of weakly open subspaces is
distributive.

\begin{lemma}
\label{lem.dist}
Let $U$, $V_1$, and $V_2$, be weakly open subspaces of a locally
noetherian space $X$. Then
$$
U \cap (V_1 \cup V_2)=(U \cap V_1) \cup (U \cap V_2).
$$
\end{lemma}
\begin{pf}
Let $W$, $Z_1$, and $Z_2$, be weakly closed subspaces of $X$ such
that $U=X \backslash W$ and $V_i=X \backslash Z_i$.
Then $V_1 \cup V_2 = X \backslash (Z_1 \cap Z_2)$, and $U \cap (V_1
\cup V_2) = X \backslash (W \bullet (Z_1 \cap Z_2))$ by
(\ref{eq.int.opens}). On the other hand $(U \cap V_1) \cup (U \cap
V_2) = X\backslash(W\bullet Z_1  \cap W \bullet Z_2)$. To prove the
result it suffices to show that $\Mod_{W\bullet Z_1  \cap W \bullet
Z_2} X= \Mod_{W \bullet (Z_1 \cap Z_2)}X$. Since $W\bullet(Z_1 \cap
Z_2) \subset W\bullet Z_1  \cap W \bullet Z_2$, one inclusion is
obvious. To prove the other it suffices to show that $\Mod W\bullet
Z_1  \cap W \bullet Z_2$ is contained in $\Mod_{W \bullet (Z_1 \cap
Z_2)}X$.

Let $M \in \Mod W\bullet Z_1  \cap W \bullet Z_2$.
There are exact sequences $0 \to L_i \to M \to N_i \to 0$, $i=1,2$,
with $L_i \in \Mod Z_i$, and $N_1,N_2 \in \Mod W$. 
The slices in the filtration $0 \subset L_1\cap L_2 \subset L_1 \subset
M$ belong to $\Mod Z_1 \cap Z_2$, $\Mod Z_1 \cap W$, and $\Mod W$
respectively, all of which are subcategories of $\Mod W \bullet
(Z_1 \cap Z_2)$, so $M \in \Mod_{W \bullet (Z_1 \cap Z_2)}X$.
This completes the proof.
\end{pf}

\begin{definition}
\label{defn.UsubW}
We say that a weakly closed subspace $W$ of $X$ {\sf contains} a 
weakly open subspace $U$ of $X$ if $\Mod U \subset \Mod W$; 
more precisely, 
if $\a:W \to X$ and $j:U \to X$ are the inclusion maps, we say that
$W$ contains $U$ if $j_*(\Mod U) \subset \a_*(\Mod W)$.
\end{definition}

\begin{lemma}
\label{lem.UsubZ}
Let $U$ and $Z$ be respectively a weakly open subspace and a weakly
closed subspace of a locally noetherian space $X$.
Let $j:U \to X$ and $\d:Z \to X$ denote the inclusion maps.
Then $U$ is contained in $Z$ if and only if there is a weak map
$\nu:U \to Z$ such that $\d\nu=j$. If that is the case, then
\begin{enumerate}
\item{}
$U$ is a weakly open subspace of $Z$;
\item{}
if $U=X \backslash W$, then $U=Z \backslash (Z \cap W)$;
\item{}
$U$ is an open subspace of $Z$ if it is an open subspace of $X$.
\end{enumerate}
\end{lemma}
\begin{pf}
If $U$ is contained in $Z$, then $\Mod U \subset \Mod Z$, so if
$\nu_*:\Mod U \to \Mod Z$ is the inclusion, then $\d_*\nu_*=j_*$.
Hence $\nu_*$ defines a weak map $\nu:U \to Z$ such that $\d\nu=j$.
Conversely, if there is a  weak map $\nu:U \to Z$ such that
$\d\nu=j$, then $\d_*\nu_* \cong j_*$, so $\d^!j_* \cong
\d^!\d_*\nu_* \cong \nu_*$ and $\d_*\d^!j_* \cong \d_*\nu_* \cong
j_*$, whence $j_*(\Mod U) \subset \d_*(\Mod Z)$; in other words $\Mod
U$ is contained in $\Mod Z$.

(1)
Because $\Mod U$ is a full subcategory of $\Mod X$ it is also a
full subcategory of $\Mod Z$. Because 
$\nu^*:=j^*\d_*$ is an exact left adjoint to $\d^!j_* \cong \nu_*$,
$U$ is a weakly open subspace of $Z$.

(2)
Since $\nu:U \to Z$ is a weakly open immersion, $\Mod U=\Mod Z/\sT$
where $\sT$ is the full subcategory of $\Mod Z$ where $\nu^*$
vanishes; however, if $M \in \Mod Z$, then $\nu^*M=0$ if and only
if $j^*M=0$, so $\sT=\Mod Z \cap \Mod_W Z$. Therefore $U=Z
\backslash (Z \cap W_{sat})$. However, by the remark after
Definition \ref{defn.sat}, $Z \backslash (Z \cap W_{sat}) =
Z \backslash (Z \cap W)$, so $U=Z \backslash (Z \cap W)$.

(3)
If $U$ is open in $X$, then $U=X \backslash W$ for some closed
subspace $W$ of $X$. Proposition \ref{prop.W.cap.Z} shows that $Z
\cap W$ is then a closed subspace of $Z$, so $Z \backslash (Z \cap
W)$ is an open subspace of $Z$.
\end{pf}

\begin{definition}
\label{defn.weak.closure}
The {\sf weak closure} of a weakly open subspace $U$ of $X$ is the 
smallest weakly closed subspace of $X$ that contains $U$; 
we denote it by $\Uol$. This makes sense because an intersection
of weakly closed subspaces is weakly closed---see the remarks after
Proposition \ref{prop.W.cap.Z}.
\end{definition}

If $U$ is a weakly open subspace of 
a locally noetherian space $X$, then $\Mod \Uol$
consists of all $X$-modules $N$ for which there exists a
$U$-module $P$ and $X$-submodules $L \subset M \subset j_*P$
such that $N \cong M/L$. It is clear that this subcategory of
$\Mod X$ contains $\Mod U$, is closed under subquotients, and
is closed under direct sums because direct sums are exact
in $\Mod X$ and $j_*$ commutes with direct sums.

\begin{lemma}
\label{lem.weak.closure}
If $U$ is a weakly open subspace of a locally noetherian 
space $X$, then 
\begin{enumerate}
\item{}
$ U \cap (X \backslash \Uol) = \phi$;
\item{}
if $U \cap V \ne \phi$ for all non-empty weakly open
subspaces $V$ in $X$, then $(\Uol)_{sat}=X$.
\end{enumerate}
\end{lemma}
\begin{pf}
(1)
Let $W$ be a saturated  weakly closed subspace of $X$ such
that $U=X \backslash W$ and write $j:U \to X$ be the inclusion map.
Let $\tau:\Mod X \to \Mod_W X$ be the $W$-torsion functor.
Let $M \in \Mod X$.
Since $M/\tau M$ embeds in $j_*j^*M$ it is a $\Uol$-module.
Hence $M \in \Mod \Uol\bullet W$, and we conclude that
$\Uol\bullet W=X$. Therefore $\phi=X \backslash(\Uol\bullet W)
= X \backslash (\Uol \cup W)=(X \backslash \Uol) \cap
(X \backslash W)$, as claimed.

(2)
The hypothesis implies that $X \backslash \Uol =\phi$, whence
$(\Uol)_{sat} =X$.
\end{pf}

If $X$ is the affine space in Example \ref{eg.bad.triang},
then $U=X \backslash q$ is an open subspace such that
$\Uol=X$, but $U \cap (X\backslash p)=U \cap q =\phi$.
It would be worthwhile to find conditions which ensure
that $U \cap V \ne \phi$ for all non-empty weakly open subspaces $V$
when $\Uol=X$.

\medskip

There is not a useful analogue of quasi-compactness for weakly open
covers: if $W_i$ is the weakly closed subspace of $\Spec
\ZZ$ defined at the beginning of section three and
$U_i$ is its weakly open complement, then $\Spec \ZZ$ is the union
of the $U_i$ but is not the union of any finite subset of these.

However, the same proof as for the commutative case shows that affine 
spaces satisfy an analogue of quasi-compactness for open covers. 

\begin{lemma}
If $\{U_i \; | \; i \in I\}$ is an open cover of an affine space
$X$, then there is a finite subset $F \subset I$ such that 
$\{U_i \; | \; i \in F\}$ is an open cover for $X$.
\end{lemma}
\begin{pf}
Suppose that $R$ is a coordinate ring for $X$, and that 
$Z_i$, $i \in I$, are closed subspaces of $X$ such that 
$U_i=X \backslash Z_i$.  Let $K_i$ be the two-sided ideal
of $R$ corresponding to $Z_i$; that is, $\Mod Z_i=\Mod R/K_i$.
Let $K$ denote the sum of all
the $K_i$. If $j_i:U_i \to X$ is the inclusion, then
$j_i^*(R/K_i)=0$, so $j_i^*(R/K)=0$. Hence $R/K=0$, and it follows
that $1$, and hence $R$, is contained in $\sum_{i \in F} K_i$ for some
finite set $F \subset I$. 

Now we show that the $U_i$, $i \in F$, provide an open cover of $X$.
If $M$ is a non-zero module such that $j_i^*M=0$ for all $i \in F$, 
then $j_i^*N=0$ for all $i \in F$ for every simple subquotient $N$ of 
$M$.  It follows that $N \in \Mod Z_i$ for all $i
\in F$, so $NK_i =0$ for all $i \in F$; this contradicts the
fact that $R=\sum_{i \in F} K_i$.
\end{pf}

\section{Containment and intersection of subspaces}
\label{sect.contains}

We have defined containment, intersection, and union,
for a pair of  weakly closed subspaces, and for a pair of
weakly open subspaces.
We would also like to define containment and intersection for pairs
of subspaces where one of the subspaces is weakly
open and the other is weakly closed.

We begin with the following question: if $Z$ and $W$ 
are weakly closed subspaces of $X$ such that $W \cap Z=\phi$, is
$Z$ a weakly closed subspace of $X\backslash W$? Unfortunately 
it need not be.

To determine exactly when $Z$ is a weakly closed subspace of
$X\backslash W$ we need to address the following general question.
If $\b:Z \to X$ is a weak map, and $j:X \backslash W \to X$
is the inclusion of a weakly open subspace, when is there a weak
map $\c:Z \to X \backslash W$ such that the diagram
$$
\begin{CD}
Z @>{\b}>> X
\\
@. @AA{j}A
\\
@.  X \backslash W
\end{CD}
$$
commutes? The case when $Z$ is a weakly closed subspace
such that $Z \cap W=\phi$ is relevant to the question posed in the
previous paragraph.

\begin{proposition}
\label{prop.contains}
Let $\b:Z \to X$ be a weak map.
Let $W$ be a weakly closed subspace of $X$; let
$\a:W \to X$ and $j:X \backslash W \to X$ be the
inclusion maps. 
The following are equivalent:
\begin{enumerate}
\item{}
$\Hom_X(\a_*M,\b_*N)=
\Ext^1_X(\a_*M,\b_*N)=0$ for all $M \in \Mod W$ and all
$N \in \Mod Z$;
\item{}
$j_*j^*\b_* \cong \b_*$;
\item{}
there is a unique weak map $\c:Z \to X \backslash W$ such that 
$\b=j\c$.
\end{enumerate}
Suppose that these conditions hold. Then
\begin{enumerate}
\item{}
if $\b$ is a map, so is $\c$;
\item{}
if $\b$ is an affine map, so is $\c$;
\item{}
if $\b$ is a weakly closed (resp. weakly open, closed) immersion,
so is $\c$.
\end{enumerate}
\end{proposition}
\begin{pf}
First we establish the uniqueness claim in part (3). If $\c$ and 
$\c'$ are weak maps from $Z$ to $X\backslash W$ such that
$\b=j\c=j\c'$, then $\b_* \cong j_*\c_* \cong j_*\c'_*$. But 
$j^*j_* \cong \id$, so $j^*\b_* \cong \c_* \cong \c'_*$. 
Because $\c_*$ and $\c'_*$ are naturally equivalent, the maps
$\c$ and $\c'$ are the same (Definition \ref{defn.map}).

Define $\c_*:\Mod Z \to \Mod X \backslash W$ by $\c_*=j^*\b_*$.
This is an exact functor so defines a weak map $\c:Z \to X
\backslash W$. (It is not necessarily true that $j\c=\b$.)
Let $\tau: \Mod X \to \Mod_W X$ denote the torsion functor that 
takes the $W$-support of 
a module. If $N \in \Mod Z$, then there is an exact sequence
\begin{equation}
\label{eq.locn}
0 \to \tau (\b_*N) \to  (\b_*N) \to j_*j^*(\b_*N) \to 
R^1\tau (\b_*N) \to 0
\end{equation}

(1) $\Rightarrow$ (2)
The vanishing of Hom implies that $\b_*N$ cannot have a non-zero
$W$-submodule, so $\tau(\b_*N)=0$. The vanishing of $\Ext^1$ 
implies that
$R^1\tau (\b_*N) =0$. Hence the natural transformation $\b_* \to
j_*j^*\b_*$ is an isomorphism.

(2) $\Rightarrow$ (3)
Since $j_*\c_* = j_*j^*\b_* \cong \b_*$, $j\c=\b$.

(3) $\Rightarrow$ (1)
If there is a weak map $\c:Z \to X \backslash W$ such that
$j\c=\b$, then $\b_* \cong j_*\c_*$, so
$j^*\b_*\cong j^*j_*\c_* \cong \c_*$, so we may assume that $\c_*=
j^*\b_*$. It follows that $j_*j^*\b_* \cong j_*\c_* \cong
\b_*$, whence (2) holds. But $j_*j^*\b_*N$ is the largest essential
extension of $\b_*N$ by a module in $\Mod_W X$. It follows that one
cannot extend $\b_*N$ in an essential way by a $W$-module. 
In other words, if $0 \to \b_*N \to E^0 \to E^1 \to \ldots$ is a
minimal injective resolution in $\Mod X$, then $\Hom_X(\a_*M,E^1)=0$
for all $W$-modules $M$; (1) follows.

Now suppose that the three equivalent conditions hold.

(1)
Suppose that $\b$ is a map. Thus $\b_*$ has a left adjoint $\b^*$.
If $L \in \Mod X \backslash W$ and $N \in\Mod Z$, then
\begin{align*}
\Hom_Z(\b^*j_*L,N) & \cong \Hom_X(j_*L,\b_*N)
\\
& \cong \Hom_X(j_*L,j_*j^*\b_*N)
\\
& \cong \Hom_{X \backslash W}(L,j^*\b_*N)
\\
& \cong \Hom_{X \backslash W}(L,\c_*N).
\end{align*}
Therefore $\c^*:=\b^*j_*$ is a left adjoint to $\c_*$,  so
$\c$ is a map.

(2)
Suppose that $\b$ is an affine map.
We have just seen that $\c$ is a map. 
Since $\b_*$ is faithful and $\b_* \cong j_*\c_*$, $\c_*$ is faithful.
Since  $\b_*$ has a right adjoint, $\b^!$ say,
$\c_*$ has a right adjoint, namely $\c^!:=\b^!j_*$.
Therefore $\c$ is an affine map.

(3)
Suppose that $\b$ is a weakly closed immersion.
Then $\b_*$ is fully faithful, $\b_*(\Mod Z)$ is closed under
subquotients in $\Mod X$, and $\b_*$ has a right adjoint $\b^!$.
Thus $\c^!:=\b^!j_*$ is right adjoint to $\c_*$.
Since $\c^!\c_*\cong \b^!j_*j^*\b_* \cong \b^!\b_*$, we have
$$
\Hom_{X \backslash  W}(\c_*L,\c_*N) \cong 
\Hom_Z(L,\c^!\c_*N) \cong \Hom_Z(L,\b^!\b_*N) = \Hom_Z(L,N).
$$
Hence $\c_*$ is fully faithful, and we can view $\Mod Z$ as a full
subcategory of $\Mod X \backslash W$ via $\c_*$.
To see that $\c_*$ is a weakly closed immersion
we must check that $\c_*(\Mod Z)$ 
is closed under subquotients in $\Mod X \backslash W$.

Let $N$ be a $Z$-module and suppose that $0 \to K \to \c_*N \to L
\to 0$ is an exact sequence in $\Mod X \backslash W$. Applying
$j_*$ produces an exact sequence
$$
0 \to j_*K \to j_*j^*\b_* N \cong \b_*N \to j_*L \to R^1j_*K
$$
Since $j_*K \to \b_*N$ is an $X$-submodule of a module in
$\b_*(\Mod Z)$, which is closed under $X$-submodules, 
$j_*K \in \b_*(\Mod Z)$; thus $j_*K \cong
\b_*\b^!j_*K$. Hence $K \cong j^*j_*K \cong j^*\b_*\b^!j_*K \cong
\c_*\c^!K$. In particular, $K \in \c_*(\Mod Z)$.
If we apply $\c_*$ to the exact sequence $0 \to \c^!K \to N \to N
/\c^!K \to 0$ we obtain an exact sequence
$0 \to \c_*\c^!K \cong K \to \c_*N \to \c_*(N/\c^!K) \to 0$. It
follows that $L \cong  \c_*(N/\c^!K) \in \c_*(\Mod Z)$.
Hence $\c_*(\Mod Z)$ is closed under subquotients.
This completes the proof that $\c$ is a weakly closed
immersion if $\b$ is.

It follows that $\c$ is a closed immersion if $\b$ is.

Finally, suppose that $\b$ is a weakly open immersion. Thus $\b_*$
has an exact left adjoint $\b^*$. Since $\b$ is a map, so is $\c$.
In other words, $\c_*$ has a left adjoint $\c^*:=\b^*j_*$.
If $M \in \Mod W$ and $N \in \Mod
Z$, then
$$
\Hom_Z(\b^*\a_*M,N) \cong \Hom_X(\a_*M,\b_*N)=0
$$
so $\b^*$ vanishes on all $W$-modules. 
Since $\b^*$ is exact and commutes with direct
limits, it therefore vanishes on $\Mod_WX$. By \cite[Coroll. 2 and
3, pp. 368-369]{G},
there is an exact functor $H:\Mod X \backslash W \to \Mod Z$
such that $Hj^* \cong \b^*$. However, $H\cong Hj^*j_* \cong
\b^*j_*=\c^*$. Therefore $\c_*$ has an
exact left adjoint, thus showing that $\c$ is a weakly open immersion.
\end{pf}

\begin{definition}
\label{defn.contains} 
Let $W$ be a weakly closed subspace of $X$. Let $V$ be either a
weakly closed or weakly open subspace of $X$. 
We say that $V$ is {\sf contained in} $X \backslash W$, and write
$V \subset X \backslash W$, if $\Mod V \subset \Mod X \backslash W$.
\end{definition}

{\bf Remarks.}
{\bf 1.}
If $\Mod V \subset \Mod X \backslash W$, then it follows from Lemma
\ref{lem.gabriel} that $\Ext_X^1(M,N)=0$ for all $W$-modules $M$ 
and all $V$-modules $N$. Hence, by Proposition \ref{prop.contains}, the
inclusion $\c_*:\Mod V \to \Mod X \backslash W$ is the direct image
functor for a weak map $\c:V \to X \backslash W$ satisfying
$\b=j\c$. Furthermore, $V$ is a weakly closed (resp., weakly open,
closed) subspace of $X \backslash W$ if it is weakly closed (resp.,
weakly open, closed) in $X$.

{\bf 2.}
Let $Z$ be a weakly closed subspace of $X$. If $Z \subset X\backslash
W$, then $Z \cap W=\phi$. However, if $Z \cap W=\phi$, then
$Z$ is contained in $X \backslash W$ if and only if 
$\Ext_X^1(M,N)=0$ for all $W$-modules $M$ and all $Z$-modules $N$.

{\bf 3.}
Example \ref{eg.bad.triang} exhibits distinct closed points 
$p,q \in X$ such that $q$ is not contained 
in $X \backslash p$. Compare this with 
Corollaries \ref{cor.pts.in.complement} and \ref{cor.pts.in.divisor}.

{\bf 4.}
If $V$ is contained in $X \backslash W$, so are the weakly closed
and weakly open subspaces of $V$.

\begin{corollary}
\label{cor.pts.in.complement}
Let $p$ be a closed point of $X$, and let $W$ be a weakly 
closed subspace of $X$ that does not contain $p$.
Let $j:X \backslash W \to X$ be the inclusion.
Suppose that $\Ext^1_X(j_*M,\cO_p)=0$ for all $W$-modules $M$. Then
\begin{enumerate}
\item{}
$p$ is contained in $X\backslash W$;
\item{}
$p$ is a closed point of $X\backslash W$;
\item{}
if $q$ is a closed point in $X\backslash W$, then
$j_*\cO_q \cong \cO_p$ if and only if $q=p$.
\end{enumerate}
\end{corollary}
\begin{pf}
Let $\b:p \to X$ denote the inclusion.
By Proposition \ref{prop.contains}, there is a closed immersion
$\c:p \to X \backslash W$ such that $\b=j\c$. Thus we can
view $p$ as a closed subspace of $X \backslash W$. 
The fact that $p$ is a closed point of $X \backslash W$ 
follows from the fact that $\c_*\cO_p=j^*\b_*\cO_p$
is simple in $\Mod X \backslash W$ (a localization functor 
sends a simple module to a simple module or zero).
Part (3) follows from the fact that $j_*$ is 
fully faithful.
\end{pf}

\begin{corollary}
\label{cor.open.in.open}
Let $W \subset Z \subset X$ be weakly closed subspaces of $X$.
Then $X \backslash Z$ is a weakly open subspace of $X \backslash W$. 
Precisely, there is a weakly open immersion $\c:X \backslash Z \to X
\backslash W$ such that the following diagram commutes:
$$
\begin{CD}
X \backslash Z @>{\b}>> X
\\
@. @AA{j}A
\\
@. X\backslash W.
\end{CD}
$$
\end{corollary}
\begin{pf}
By Lemma \ref{lem.gabriel}, $\Ext^1_X(-,\b_*N)$ vanishes on all
$Z$-modules. 
Thus the first condition in Proposition \ref{prop.contains} is
satisfied, and it follows that such a $\c$ exists.
\end{pf}

{\bf Remarks. 1.}
By Proposition \ref{prop.int.opens},
if $U$ and $V$ are weakly open subspaces of $X$, then
there is a commutative diagram of weakly open immersions
\begin{equation}
\label{eq.cap.opens}
\begin{CD}
U \cap V @>{\a}>> U
\\
@V{\b}VV @VV{i}V
\\
V@>>{j}> X.
\end{CD}
\end{equation}
Corollary \ref{cor.open.in.open} also confirms this:
let $W$ and $Z$ be weakly closed subspaces such that
$U=X \backslash W$ and $V=X \backslash Z$;
let $\c:U \cap V = X \backslash (W \cup Z) \to X$ be the inclusion;
since $W \subset W \cup Z$, there is a weakly 
open immersion $\a:X \backslash (W \cup Z) \to X \backslash W$ 
such that $i\a=\c$ (by the proof of Proposition
\ref{prop.contains}, $\a$ is obtained by defining $\a_*=i^*\c_*$);
similarly, there is a weakly open immersion
$\b:X \backslash (W \cup Z) \to X \backslash Z$ such that
$j\b=\c$.

{\bf 2.}
In contrast to the case for schemes \cite[Ch. II, Ex. 4.3]{H},
we cannot conclude that $U \cap V$ is affine even if $U$ and $V$ are
open affine subspaces of an affine space $X$. 
This is not a surprise. Suppose that
$X$ is affine with coordinate ring $R$. If $i$ and $j$ are affine
maps, then $U$ and $V$ are affine, with coordinate rings $R_U$ and
$R_V$ say; we can choose $R$ such that $i$ and $j$ are induced by
ring homomorphisms $R \to R_U$ and $R \to R_V$, both of which are
necessarily epimorphisms in the category of rings. However, when
$R$ is not commutative, there may not be a ring structure on $R_U
\otimes_R R_V$ such that the natural map $R \to R_U \otimes_R R_V$ 
is a ring homomorphism.

\begin{proposition}
\label{prop.fiber}
Let $U$ and $V$ be weakly open subspaces of a locally noetherian
space $X$. Then $U \cap V$ is the fiber product $U \times_X V$ of
$U$ and $V$ in the category of spaces (with maps as the morphisms).
\end{proposition}
\begin{pf}
Consider the diagram (\ref{eq.cap.opens}) above.
Suppose that $f:Y \to U$ and $g:Y \to V$ are maps such that
$if=jg$. We must show that there is a unique map $h:Y \to U \cap V$
such that $f=\a h$ and $g=\b h$.

The inclusion $\c:U \cap V \to X$ is given by
$\c_*=i_*\a_*=j_*\b_*$.
Since $\c^* \cong \a^*i^*$, we have $\c^*i_* \cong \a^*$.

Because the image of $i_*f_*$ is the same as that of $j_*g_*$, 
it is contained in $\Mod
U \cap \Mod V = \Mod U \cap V$, whence $i_*f_* \cong
\c_*\c^*i_*f_*$. Thus
$$
f_* \cong i^*i_*f_* \cong i^*\c_*\c^*i_*f_* \cong
i^*(i_*\a_*\a^*i^*)i_*f_* \cong \a_*\a^*f_*.
$$
By Proposition \ref{prop.contains}, there is a unique map 
$h:Y \to U \cap V$ such that $f=\a h$.
The proof of that result shows that $h_*=\a^*f_*$.
The ``same'' argument shows there is a unique map $h':Y \to U \cap
V$ such that $g=\b h'$, and $h'_*=\b^*g_*$. 
Therefore
$$
h^*=\a^*f_* \cong \c^*i_*f_* \cong \c^*j_*g_* \cong \b^*g_* =h'_*,
$$
and we conclude that $h=h'$, completing the proof.
\end{pf}

We now generalize the question posed prior to Proposition
\ref{prop.contains}.
If $Z$ and $W$ are weakly closed subspaces of a space $X$, when is
$Z \backslash (Z \cap W)$ a weakly closed subspace of $X \backslash W$?
Example \ref{eg.bad.triang} shows that the answer is not ``always''.

\begin{proposition}
\label{prop.intersect}
Let $W$ and $Z$ be weakly closed subspaces of a locally noetherian space
$X$. Denote the inclusion maps by
$$
\begin{CD}
Z @>{\d}>> X @<{\a}<< W.
\\
@A{\ve}AA @AA{j}A
\\
Z \backslash (Z \cap W) @. X\backslash W
\end{CD}
$$
There is a weakly closed immersion 
$\c:Z \backslash (Z \cap W) \to X\backslash W$ such that $j\c=\d\ve$
if and only if $\Ext^1_X(\a_*M,\d_*\ve_*N)=0$ for all $M \in \Mod W$
and all $N \in \Mod Z \backslash (Z \cap W)$.
In that case, if $\d$ is a closed immersion, so is $\c$.
\end{proposition}
\begin{pf}
We set $\b=\d\ve$. 
Thus $\b:Z \backslash(Z\cap W) \to X$ is a weak map.

($\Rightarrow$)
Suppose there is a weak map $\c$ such that $j\c=\d\ve$.
By Proposition \ref{prop.contains}, $\Ext^1_X(\a_*M,\d_*\ve_*N)=0$
for all $M$ and $N$.

($\Leftarrow$)
Proposition \ref{prop.contains} immediately yields
a weak map $\c:Z \backslash (Z \cap W)
\to X\backslash W$ such that $j\c=\b=\d\ve$. Explicitly, $\c_*=j^*\b_*=
j^*\d_*\ve_*$. Furthermore, $\b_* \cong j_*j^*\b_*$.

To show that $\c$ is a weakly closed immersion we must show that
that $\c_*(\Mod Z \backslash (Z \cap W))$
is closed under subquotients in $\Mod X\backslash W$ and 
that $\c_*$ is exact, fully faithful, and has a right adjoint.

Since it is a composition of left exact functors, $\c_*$ is left
exact. To show that $\c_*$ is exact, it suffices to show that
$j^*\d_*(R^1\ve_*)=0$. 
Now $R^1\ve_*$ takes values in $\Mod_{Z \cap W} Z$ and
$j^*\d_*$ vanishes on this because $j^*$ vanishes on $\Mod_W X$.
Hence $j^*\d_*(R^1\ve_*)=0$ and  $\c_*$ is exact.

Because $j^*$ and $\d_*$ have right adjoints, they commute with
direct limits. By Lemma \ref{lem.gabriel}, $\ve_*$ commutes with
direct limits. Therefore $\c_*$ commutes with direct limits. Since $\c_*$
is exact it therefore has a right adjoint, which we denote by
$\c^!$.

It is useful to have a precise formula for $\c^!$.
Since $j^*\d_*$ vanishes on $\Mod_{Z \cap W} Z$, 
$j^*\d_*$ vanishes on the kernel and cokernel of
the natural transformation $\id_Z \to \ve_*\ve^*$. Hence the natural
transformation $j^*\d_* \to j^*\d_*\ve_*\ve^*$ is an isomorphism.
That is, $j^*\d_* \cong \c_*\ve^*$. Taking right adjoints, $\d^!j_*
\cong \ve_*\c^!$. It follows that 
$$
\c^! \cong \ve^*\d^!j_*.
$$

The fact that $\c_*$ is fully faithful follows from the
calculation
\begin{align*}
\Hom_{Z \backslash (Z \cap W)}(L,N) 
& \cong 
\Hom_{Z \backslash (Z \cap W)}(\ve^*\ve_*L,N) 
\\
& \cong 
\Hom_{Z}(\ve_*L,\ve_*N) 
\\
& \cong 
\Hom_{Z}(\ve_*L,\d^!\d_*\ve_*N) 
\\
& \cong 
\Hom_{X}(\d_*\ve_*L,\b_*N) 
\\
& \cong 
\Hom_{X}(\d_*\ve_*L,j_*j^*\b_*N) 
\\
& \cong 
\Hom_{X \backslash W}(j^*\d_*\ve_*L,j^*\b_*N) .
\end{align*}
Because $\c_*$ is fully faithful, the unit $\id_{Z \backslash (Z
\cap W)} \to \c^!\c_*$ is an isomorphism.

Now we check that $\c_*(\Mod Z \backslash (Z \cap W))$ is closed
under subquotients in $\Mod X\backslash W$. 
Let $N \in \Mod Z \backslash (Z
\cap W)$, and consider an exact sequence $0 \to P \to \c_*N \to Q
\to 0$ in $\Mod X\backslash W$. Hence
$j_*P$ is an $X$-submodule of $j_*\c_*N \cong \d_*\ve_* N$. But
$\d_*(\Mod Z)$ is closed under submodules, so $j_*P$ is a
$Z$-module. Formally, $j_*P \cong \d_*\d^!j_*P$, whence
$$
P \cong j^*j_*P \cong j^*\d_*\d^!j_*P \cong \c_*\ve^*\d^!j_*P \cong
\c_*\c^! P.
$$
Thus $P$ is in the image of $\c_*$.
Now, if we apply $\c_*$ to the exact sequence $0 \to \c^!P \to
\c^!\c_* N \cong N \to N/\c^!P \to 0$, we obtain the original
sequence $0 \to \c_*\c^!P \cong P \to \c_*N \to \c_*(N/\c^!P) \to
0$, from which it follows that $Q \cong \c_*(N/\c^!P)$. 

Finally, we show that if $\d$ is a closed immersion, so is $\c$. If
$\d^*$ is a left adjoint to $\d_*$, then
\begin{align*}
\Hom_Z(\ve^*\d^*j_* P, L) 
& 
\cong \Hom_X(j_*P,\d_*\ve_*L) \cong \Hom_X(j_*P,j_*j^*\b_*L) 
\\
&
\cong \Hom_{X\backslash W}(j^*j_*P,j^*\b_*L) \cong \Hom_{X\backslash
W}(P,\c_*L).
\end{align*}
Hence $\c^*:=\ve^*\d^*j_*$ is a left adjoint to $\c_*$, 
showing that $\c$ is a closed immersion.
\end{pf}

\begin{lemma}
\label{lem.open.cap.closed}
In the setting of Proposition \ref{prop.intersect}
$$
j_*\c_*(\Mod Z \backslash (Z \cap W)) = \d_*(\Mod Z) \cap j_*(\Mod
X \backslash W).
$$
\end{lemma}
\begin{pf}
Since $\d_*$ and $j_*$ are fully faithful, the category on the
right-hand side is a full subcategory of $\Mod X$. So is the
category on the left-hand side,
 because $\c_*$ is fully faithful. Hence to
verify this equality it is enough to check it on objects. Because
$j_*\c_*\cong \d_*\ve_*$, the left-hand side is contained in the
right-hand side. Conversely, a module on the right-hand side is of
the form $\d_*L$ for some $L \in \Mod Z$ and it satisfies $\d_*L
\cong j_*j^*\d_*L$. The proof of Proposition \ref{prop.intersect}
showed that $j^*\d_* \cong \c_*\ve^*$, so
$\d_*L \cong j_*j^*\d_*L \cong j_*\c_*\ve^*L$, which is in the
left-hand side.
\end{pf}

\begin{definition}
\label{defn.Z.cap.U}
Let $W$ and $Z$ be weakly closed subspaces of a 
locally noetherian space $X$, and set $U=X \backslash W$. 
Let $\a:W \to X$ and $\b:Z \backslash (Z \cap W) \to Z \to X$ be
the inclusions.
If $\Ext^1_X(\a_*M,\b_*N)=0$ for all $M \in \Mod W$ and all
$N$ in $\Mod Z \backslash (Z \cap W)$
we define $Z \cap U$ by declaring
$$
\Mod Z \cap U := \Mod Z \cap \Mod U.
$$
\end{definition}

{\bf Remarks.}
In these remarks, $W$ and $Z$ are weakly closed subspaces of 
$X$, $U$ denotes $X \backslash W$, and we use the notation of
Proposition \ref{prop.intersect}.

{\bf 1.}
Example \ref{eg.bad.triang} shows that $Z \cap U$ is not always
defined. In that example, $Z \cap (X \backslash q)$ is not defined
because $Z \backslash (Z \cap q)=p$ and $\Ext^1_X(\cO_q,\cO_p) \ne
0$. This shows that $Z \cap U$ need not be defined even if $U \to
X$ is affine.

{\bf 2.}
If $Z$ is contained in $U$ (Definition \ref{defn.contains}), then
Lemma \ref{lem.gabriel} ensures that $Z \cap U$ is defined, and $Z
\cap U = Z$. 

{\bf 3.}
If $U$ is contained in $Z$ (Definition \ref{defn.UsubW}), then
$Z \backslash (Z \cap W)=U$ by Lemma \ref{lem.UsubZ}, so
Lemma \ref{lem.gabriel} ensures that $Z \cap U$ is defined, 
and $Z \cap U = U$. In particular, $X \cap U=U$ and $\Uol \cap U
=U$. 

{\bf 4.}
Lemma \ref{lem.gabriel} implies that $W \cap U$ is defined, and $W
\cap U=\emptyset$.

{\bf 5.}
If $V \subset U$ is a weakly open subspace of $U$, then $V$ is
weakly open in $X$, and $\Vol \subset \Uol$ where these are the
weak closures in $X$. We do not know if $\Vol \cap U$ is defined.

{\bf 6.}
If $U$ is also weakly closed in $X$, then we may use Definition
\ref{defn.intersect} to define $Z \cap U$; it is a weakly closed
subspace of $Z$, of $U$, and of $X$, and is also weakly open in $Z$
by Proposition \ref{prop.W.cap.Z}. We claim that $Z \cap U$ is
also defined according to Definition \ref{defn.Z.cap.U}; comparing
Definitions \ref{defn.intersect} and \ref{defn.Z.cap.U} it is then
clear that the two definitions of $Z \cap U$ agree. To verify the
claim we must show, in the notation of  Proposition
\ref{prop.intersect}, that there is a map $\c:Z\backslash(Z \cap W)
\to X \backslash W$ making the diagram commute. By Proposition
\ref{prop.contains}, it suffices to show that $j_*j^*\d_*\ve_* \cong
\d_*\ve_*$. Since $j:U \to X$ is a weakly closed immersion, the
natural map $\d_*\ve_*N \to j_*j^*\d_*\ve_*N$ is epic.
The kernel is $\tau(\d_*\ve_*N)$, where $\tau$ is the $W$-torsion
functor. However, zero is the only $Z \cap W$-submodule of
$\ve_*N$, so zero is the only $W$-submodule of $\d_*\ve_* N$,
whence $\tau(\d_*\ve_*N)=0$.

{\bf 7.}
If $Z$ is also weakly open in $X$, then we may use Definition
\ref{defn.intersect.opens} to define $Z \cap U$; we claim that 
$Z \cap U$ is
also defined according to Definition \ref{defn.Z.cap.U}; comparing
Definition \ref{defn.Z.cap.U} with Proposition
\ref{prop.int.opens}, it is then
clear that the two definitions of $Z \cap U$ agree. To verify the
claim we must show that $\Ext^1_X(\a_*M,\d_*\ve_*N)=0$ for all $M
\in \Mod W$ and all $N \in \Mod Z \backslash (Z \cap W)$. Since $Z$
is both weakly open and weakly closed, $\d_*$ is exact. Therefore
the spectral sequence (\ref{eq.spec.seq.open}) for $\d:Z \to X$
collapses, and gives $\Ext^1_X(\a_*M,\d_*\ve_*N) \cong
\Ext^1_Z(\d^*\a_*M,\ve_*N)$. Since $\d^*\a_*M$ is a $Z \cap
W$-module, Lemma \ref{lem.gabriel} applied to the open subspace $ Z
\backslash (Z \cap W)$ of $Z$ shows that this $\Ext^1_Z$ group
vanishes.

\begin{proposition}
\label{prop.Z.cap.U}
Let $W$ and $Z$ be weakly closed subspaces of a 
locally noetherian space $X$, and set $U=X \backslash W$. 
Suppose that $Z \cap U$ is defined. Then
\begin{enumerate}
\item{}
$Z \cap U = Z \backslash (Z \cap W)$;
\item{}
$Z \cap U$ is a weakly open subspace of $Z$ and
a weakly closed subspace of $U$;
\item{}
if $Z$ is closed in $X$, then $Z \cap U$ is closed in $U$;
\item{}
if $U \to X$ is an affine map, so is $Z \cap U \to Z$.
\end{enumerate}
\end{proposition}
\begin{pf}
Part (1) follows from Lemma \ref{lem.open.cap.closed}, and parts
(2) and (3) follow from Proposition \ref{prop.intersect}.

(4)
In the notation of Proposition \ref{prop.intersect} we must show
that $\ve_*$ is faithful and has a right adjoint, where $\ve$ is
the weakly open immersion $Z \cap U \to Z$. We already know that
$\ve_*$ is faithful and commutes with direct sums since $\ve^*$ is
a localization. Since $j:U \to X$ is affine, $j_*$ and hence
$j_*\c_*$ is exact. But $j_*\c_* \cong \d_*\ve_*$; since a sequence
of $Z$-modules is exact if and only if it is exact when considered
as a sequence of $X$-modules, it follows that $\ve_*$ is exact.
\end{pf}

The language we have been using is close to that of topology. 
For this paragraph only we will say that the notions of weakly
closed and weakly open subspaces provide a space with an
almost-topological structure. If $Z$ and $U$ are respectively a
weakly closed and weakly open subspace of $X$, it (almost) makes
sense to ask if the almost-topological structures on $Z$ and $U$
(which they have by virtue of the fact that they are spaces) are
the same as the induced  almost-topological structures on them.
The next result shows this is true for $U$, but the example that
follows it shows this is false for $Z$.

\begin{proposition}
\label{prop.induced.top2}
Let $U$ be a weakly open subspace of a locally noetherian space $X$. 
Then every weakly closed subspace of $U$ is of the form 
$Z \cap U$ for some weakly closed subspace $Z$ of $X$.
\end{proposition}
\begin{pf}
Let $V$ be a weakly closed subspace of $U$, and let $\c:V \to U$
and $j:U \to X$ be the inclusions. We seek a weakly closed subspace
$Z$ of $X$ such that $V=Z \cap U$. That is, if $\d:Z \to X$ denotes
the inclusion, we seek a commutative diagram
$$
\begin{CD}
Z\cap U=V @>{\c}>> U
\\
@V{\ve}VV @VV{j}V
\\
Z@>>{\d}> X.
\end{CD}
$$
By hypothesis, $j_*$ and $\c_*$ are fully faithful, so we can think
of them as inclusions of full subcategories.

Let $\sA$ be the full subcategory of $\Mod X$ consisting of all
$X$-modules $N$ that can be written as $N \cong M/L$ where $L
\subset M \subset j_*\c_* P$ are $X$-submodules of $j_*\c_*P$ for
some $P \in \Mod V$. 

Clearly $\sA$ is closed under subquotients in $\Mod X$.
Hence the inclusion $\d_*:\sA \to \Mod X$ is exact. If 
$N_i=M_i/L_i$, $i \in I$, belong to $\sA$ where $L_i \subset M_i
\subset j_*\c_*P_i$, then taking direct sums in $\Mod X$, we have
$\oplus_i N_i \cong \oplus_i M_i/\oplus_i L_i$, and 
$$
\bigoplus_i L_i \subset \bigoplus_i M_i \subset 
\bigoplus j_*\c_* P_i \cong
j_*\c_*(\bigoplus_i P_i)
$$ 
where the final direct sum is taken in $\Mod V$ and the final
isomorphism is due to the fact that $j_*$ and $\c_*$ commute with
direct limits. It follows that $\oplus_i N_i$ is in $\sA$.
Thus $\sA$ is closed under direct limits in $\Mod X$. Hence 
$\d_*$ commutes with direct limits, and therefore
has a right adjoint $\d^!$. 
So we can define a weakly closed subspace $Z$ of $X$ by declaring
$\Mod Z$ to be $\sA$.

We now show that $j_*\c_*(\Mod V)=\d_*(\Mod Z) \cap j_*(\Mod U)$.
One inclusion is easy: if $P \in \Mod
V$, then $j_*\c_*P$ is in $\Mod Z$ by definition, thus showing that
the right-hand side contains the left-hand side. Now suppose that
$N$ is in the right-hand side. Thus $N=M/L$ where $L \subset M
\subset j_*\c_*P$ for some $P \in \Mod V$. 
There are inclusions $j^*L \subset j^*M \subset j^*j_*\c_*P \cong
\c_*P$. Since $\c_*(\Mod V)$ is closed under subquotients in $\Mod
U$, $j^*M/j^*L$ is in $\c_*(\Mod V)$. Hence $j^*N \in \c_*(\Mod V)$.
But $N$ is in $j_*(\Mod U)$, so $N \cong j_*j^*N \in j_*\c_*(\Mod
V)$. 

To show that $Z \cap U$ is defined we must show that
$\Ext^1_X(\a_*-,\d_*\ve_*-)=0$. Since $\d_*\ve_*=j_*\d_*$ this
follows from Lemma \ref{lem.gabriel}. Hence $Z \cap U$ is defined,
and it follows from the previous paragraph and Lemma
\ref{lem.open.cap.closed} that $V=Z \cap U$.
\end{pf}

We should probably call the $Z$ constructed in the previous 
proof the {\sf weak closure} of $V$ in $X$, and extend Definition
6.19 accordingly.  If we do that, then Proposition 
\ref{prop.induced.top2}
says that $V=\Vol \cap U$ whenever $V$ is a weakly closed subspace
of a weakly open subspace $U$. Thus $V$ is the analogue of 
a locally closed subscheme in the sense that it is a weakly 
open subspace of its weak closure $\Vol$.

If products are exact in $\Mod X$ (for example, if $X$ is affine),
there is a version of Proposition \ref{prop.induced.top2} with
``weakly closed'' replaced by ``closed''. In that case, the
equality $\Vol \cap U = V$ extends the result that if
$R_{\cS}$ is an Ore localization of a ring $R$, and $J$ is a
right ideal of $R_{\cS}$, then $J=(J \cap R)R_{\cS}$.

\begin{example}
\label{eg.bad.triang2}
We exhibit a closed subspace $Z \subset X$ and an open
subspace of $Z$ that is not of the form $Z \cap U$ for any weakly open
subspace $U$ of $X$.
We retain the notation in Example \ref{eg.bad.triang}.
The weakly open subspaces of $X$ are 
$$
\phi, \, X\backslash p=q, \, X\backslash q, \, X
$$
and the weakly open subspaces of $Z$ are 
$$
\phi, \, Z \backslash p =q, \, Z \backslash q =p, \, Z.
$$
Therefore  $p$, which is weakly open (even open) in $Z$, is not 
equal to $Z \cap U$ for any weakly open subspace $U$ in $X$.
\end{example}

\section{Effective Divisors}
\label{sect.div}

Van den Bergh \cite{vdB} has defined the notion of a divisor on a
non-commutative space. The reader is referred there for the
definition and basic results. Further properties of effective 
divisors can be found in \cite{J}.

\begin{definition}
\cite{J}
\label{defn.div}
An {\sf effective divisor} $Y$ on a non-commutative space $X$ 
is an invertible proper subobject $o_X(-Y)$ of $o_X$ in the category of
$X$-$X$-bimodules.
\end{definition}

If $X$ is a scheme such that $\Qcoh X$ is a Grothendieck category, 
then there is a bijection between effective Cartier
divisors and effective divisors in the sense
of Definition \ref{defn.div}.

If $Y$ is an effective divisor, then $o_X(-Y)$ is a proper ideal 
in $o_X$, so the results in \cite[Section 3.5]{vdB} show that
$Y$ determines, and is determined by, a non-empty closed 
subspace of $X$ which we also denote by $Y$. The category $\Mod Y$
is the full subcategory of $\Mod X$ consisting of those $M$ such
that the natural map $M(-Y) \to M$ is zero.

We write $i:Y \to X$ for the inclusion of the effective divisor.
The left and right adjoints to the inclusion $i_*:\Mod Y \to \Mod X$
are related via their derived functors as follows:
\begin{equation}
\label{eq.derived.functors}
R^1i^! \cong i^*(-)(Y)
\qquad \hbox{and} \qquad
L_1i^* \cong i^!(-)(-Y).
\end{equation}
Because $i^*$ is right exact, $R^2i^!=0$. Therefore
the spectral sequence (\ref{eq.sp.seq}) degenerates,
giving a long exact sequence
\begin{align}
\label{eq.five.seq}
0 \to & \Ext^1_Y(M,i^{\cstar}N)  \to \Ext^1_X(i_*M,N) \to
\Hom_Y(M,R^1i^{\cstar}N) \to
\\
\to &  \Ext^2_Y(M,i^{\cstar}N) \to \Ext^2_X(i_*M,N)
\to \Ext^1_Y(M,R^1i^{\cstar}N) \to \cdots
\notag
\end{align}

The associated exact sequence of bimodules is
\begin{equation}
\label{eq.divisor.seq}
\begin{CD}
0 @>>> o_X(-Y) @>>> o_X @>{\varphi}>> o_Y @>>> 0,
\end{CD}
\end{equation}
where $-\otimes_X o_Y=i_*i^*$ and $\cHom_X(o_Y,-)=i_*i^!$.
Applying $M \otimes_X -$ to this produces an exact sequence 
\begin{equation}
\label{eq.V.crit}
0 \to \cTor^X_1(M,o_Y)=i^!M(-Y) \to M(-Y) @>{\psi}>> M \to i^*M \to 0.
\end{equation}
We note that $\cTor^X_r(-,o_Y)=0$ for $r \ge 2$.

\begin{proposition}
Let $i:Y \to X$ and $\a:W \to X$ be the inclusions of  an effective
divisor and a weakly closed subspace respectively.
Suppose that $M(Y)$ and $M(-Y)$ belong to $\Mod W$ whenever $M$
does.
Then $W \cap Y$ is an effective divisor on $W$ if and only if 
$W \cap Y \ne \phi$ and $\cTor^X_1(o_Y,o_W)=0$.
The last condition is equivalent to 
$i^*\a_*\a^!$ vanishing on all injective $X$-modules.
\end{proposition}
\begin{pf}
By Proposition \ref{prop.W.cap.Z}, 
there is a commutative diagram 
$$
\begin{CD}
W \cap Y @>\c>> Y
\\
@V{\d}VV @VV{i}V
\\
W @>>{\a}> X.
\end{CD}
$$
in which $W \cap Y$ is a closed subspace of $W$. Furthermore,
$\c_*\d^* \cong i^*\a_*$ and $\d_*\c^! \cong \a^!i_*$.

Before proceeding we establish the equivalence claimed 
in the last sentence of the Proposition.
By definition of $\cTor^X$, if $E$ is an injective $X$-module, then
$$
\cHom_X(\cTor^X_1(o_Y,o_W),E) = \cExt^1_X(o_Y,\cHom_X (o_W,E))
=R^1i^!(\a_*\a^!E).
$$
Therefore $\cTor^X_1(o_Y,o_W)=0$ if an only if $i^*\a_*\a^!E=0$ for all
injective $X$-modules $E$. 

Now suppose that $W \cap Y$ is an effective divisor on $W$. Since
an effective divisor is non-empty, $W \cap Y \ne \phi$.
Furthermore,
$\d^* \cong R^1\d^!(-)(W \cap Y)$, so $\d^*$ vanishes on all
injective $W$-modules. Since $\a^!$ sends injectives to injectives,
$\c_*\d^*\a^!$ vanishes on all injective $X$-modules. 
Hence $i^*\a_*\a^!$ vanishes on all injective $X$-modules. 

Conversely, suppose that $W \cap Y \ne \phi$ and that
$\cTor^X_1(o_Y,o_W) =0$. Because $W \cap Y \ne \phi$, $o_{W \cap Y}
\ne 0$. Hence to show that $W \cap Y$ is an effective divisor on 
$W$, we must
prove the invertibility of the ideal $I$ of $o_W$ that is 
the kernel in the exact
sequence $0 \to I \to o_W \to o_{W \cap Y} \to 0$ of
$W$-$W$-bimodules. 
By definition, $\cHom_W(o_{W \cap Y},M)=\d_*\d^!M$ for all $M \in
\Mod W$.
The ideal $I$ is defined by the requirement that if $E$ is an injective
$W$-module, then the sequence
$$
0 \to \d_*\d^! E \to E \to \cHom_W(I,E) \to 0
$$
is exact. Equivalently, $\cHom_W(I,E)$ is defined by the
requirement that 
$$
0 \to \a_*\d_*\d^! E \to \a_*E \to \a_*\cHom_W(I,E) \to 0
$$
is exact for all injective $W$-modules $E$. But $\a_*\d_*\d^! E
\cong \cHom_X(o_Y,\a_*E)$, so $\a_*\cHom_W(I,E)$
is the cokernel of the map 
$\cHom_X(o_Y,\a_*E) \to \cHom_X(o_X,\a_*E)$. However, the
exact sequence
$$
0 \to o_X(-Y) \to o_X \to o_Y \to 0
$$
gives rise to the exact sequence
\begin{align*}
0 \to & \cHom_X(o_Y,\a_*E) \to \cHom_X(o_X,\a_*E) \to
\cHom_X(o_X(-Y),\a_*E) \to 
\\
& \cExt^1_X(o_Y,\a_*E) \to 0.
\end{align*}
Therefore $\a_*\cHom_W(I,E)$ is characterized by the fact that 
it is the kernel in the exact sequence
$$
0 \to \a_*\cHom_W(I,E) \to  \cHom_X(o_X(-Y),\a_*E) \to
\cExt^1_X(o_Y,\a_*E) \to 0.
$$
However, if $E'$ is an injective envelope of $\a_*E$ in $\Mod X$,
then $E \cong \a^!E'$, so
$$
\cExt^1_X(o_Y,\a_*E) \cong \cExt^1_X(o_Y,\cHom_X(o_W,E')) =
\cHom_X(\cTor_1^X(o_Y,o_W),E')
$$
which is zero by hypothesis.
It follows that 
$$
\a_*\cHom_W(I,E) \cong \cHom_X(o_X(-Y),\a_*E)
$$
and hence that 
$$
\cHom_W(I,-) \cong \a^! \cHom_X(o_X(-Y),-) \a_*
$$
on injective $W$-modules, and hence on all $W$-modules.

However,
if $F,G \in \Lex(X,X)$ are mutual quasi-inverses and each sends
$\Mod W$ to itself, then $\a^!F\a_*$ and $\a^!G\a_*$ are in
$\Lex(W,W)$ and are mutual quasi-inverses. 
Therefore $\cHom_W(I,-)$ is an auto-equivalence of $\Mod W$.
%
\end{pf}

The next two lemmas will be used later.

\begin{lemma}
\label{lem.VV}
\cite{VV}
If $Y$ is an effective divisor on $X$, then $\Mod_YX$ is closed under
injective envelopes in $\Mod X$.
\end{lemma}
\begin{pf}
This is proved in \cite[Proposition 8.4]{VV}. In their notation,
$-\otimes_X o_X(-Y)$ corresponds to the functor $G$, 
and $\Mod Y$ corresponds to $\cB=\cD$.
\end{pf}

\begin{lemma}
\label{lem.Ri.omega.tau}
Let $Z$ be a weakly closed subspace of a locally noetherian space
$X$. Let $j:X \backslash Z \to X$ be the inclusion map. 
Let $\tau:\Mod X \to \Mod_Z X \to \Mod X$ be the functor that takes
$Z$-torsion.
If $\Mod_Z X$ is closed under injective envelopes in $\Mod X$, then 
\begin{enumerate}
\item{}
every injective $X$-module is a direct sum of a torsion injective and
a torsion-free injective;
\item{}
for $i \ge 1$, the right-derived functors satisfy 
$
R^{i+1}\tau M \cong R^{i}j_*(j^* M).
$
\end{enumerate}
\end{lemma}
\begin{pf}
(1)
Let $E$ be an injective $X$-module. Since $E$ contains a copy of the
injective envelope of $\tau E$, and since that injective is torsion
by hypothesis, $\tau E$ is injective. Hence $E \cong \tau E \oplus Q$
with $Q$ a torsion-free injective.

(2)
Let $M \to E^{\smallbullet}$ be an injective resolution of $M$. For
each $i$, write $I^i=\tau(E^i)$, and set
$Q^i=E^i/I^i$. Then there is an exact sequence of complexes
$$
0 \to I^{\smallbullet} \to E^{\smallbullet} \to Q^{\smallbullet}
\to 0
$$
which gives a long exact sequence 
in homology.  However, $h^i(I^{\smallbullet})=R^i\tau M$, and
$h^i(E^{\smallbullet})=0$ for $i \ge 1$, so
$R^{i+1}\tau M \cong h^i(Q^{\smallbullet})$ for $i \ge 1$.

Since $Q^i$ is torsion-free, $j_* Q^i$ is an injective $X
\backslash Z$-module, and $j_*j^* Q^i \cong Q^i$. 
Since $j^*$ is exact, $j^* M \to j^* E^{\smallbullet}$
is an injective resolution in $\Mod X \backslash Z$.
The complexes $j^* Q^{\smallbullet}$ and $j^* E^{\smallbullet}$ 
are isomorphic, so $j^* M \to j^* Q^{\smallbullet}$ is an injective
resolution in $\Mod X \backslash Z$. Thus
$
R^ij_*(j^* M) \cong h^i(j_* j^* Q^{\smallbullet}) \cong
h^i(Q^{\smallbullet}).
$
\end{pf}

\begin{proposition}
\label{prop.open}
Let $Y$ be an effective divisor on a locally noetherian space $X$.
Let $i:Y \to X$ and $j:X \backslash Y \to X$ be the inclusions. 
Then 
\begin{enumerate}
\item{}
$j$ is an affine map;
\item{}
$i^!j_*=i^*j_*=0$;
\item{}
whenever $M \in \Mod X$  and $N \in \Mod X \backslash Y$,
\begin{equation}
\label{eq.div.compl}
\Ext^n_{X\backslash Y}(j^*M,N) \cong \Ext^n_X(M,j_*N).
\end{equation}
\end{enumerate}
\end{proposition}
\begin{pf}
Let $\tau:\Mod X \to \Mod_Y X \to \Mod X$ denote the 
$Y$-torsion functor.

(1)
By Lemma \ref{lem.Ri.omega.tau}, $(R^1j_*)j^* \cong R^2\tau$.
But $\tau$ is the direct limit of the functors 
$\cHom_X(o_X/o_X(-nY),-)$ for $n \ge 0$, so
$R^2\tau$ is the direct limit of the functors 
$\cExt_X^2(o_X/o_X(-nY),-)$, each of which 
is zero because $o_X(-Y)$ is an invertible bimodule. Therefore $0= 
R^2\tau\cong (R^1j_*)j^*$. Since every $X \backslash Y$-module is
of the form $j^*M$, it follows that $R^1j_*=0$ and we conclude that
$j_*$ is exact. But $j_*$
commutes with direct sums (Proposition \ref{prop.gabriel}), so it
has a right adjoint. Since $j^*j_* \cong \id_{X \backslash Y}$,
$j_*$ is fully faithful. Hence $j$ is an affine map.

(2)
The definition of $j_*$ ensures that $i^!j_*$ is zero.
Let $N \in \Mod X \backslash Y$, and suppose that $i^*j_*N$ is non-zero.
By (\ref{eq.V.crit}), there is an exact sequence
$$
0 \to (j_*N)(-Y) \to j_*N \to i^*j_*N \to 0.
$$
It follows from the construction of $j_*$, that $j_*N$ is not an
essential extension of $(j_*N)(-Y)$. Hence, if $V$ denotes the image of
$(j_*N)(-Y)$ in $j_*N$, $V \cap L=0$ for some
non-zero submodule $L \subset j_*N$
. But then
$L$ embeds in $i^*j_*N$, so is a $Y$-module,
contradicting the fact that $\tau(j_*N)=0$. From this contradiction
we deduce that $i^*j_*N=0$.

(3)
Because $j_*$ is exact, the 
spectral sequence (\ref{eq.spec.seq.open}) for $U=X \backslash Y$
collapses, giving $\Ext^n_X(-,j_*N) \cong 
\Ext^{n}_{X \backslash Y}(j^*-,N)$.
\end{pf}

\begin{proposition}
\label{prop.contained.in.divisor}
Let $Y$ be an effective divisor on $X$.
If $Z$ is a weakly closed subspace of $X$ 
such that $Z \cap Y =\phi$, then 
\begin{enumerate}
\item{}
$Z$ is contained in $X\backslash Y$, and
\item{}
$Y$ is contained in $X\backslash Z$.
\end{enumerate}
\end{proposition}
\begin{pf}
Let $i:Y \to X$ and $\b:Z \to X$ denote the inclusions.
Because $Z \cap Y=0$, both $i^!$ and $i^*$ vanish on all
$Z$-modules. 
Let $M$ and $N$ be respectively a $Y$-module and $Z$-module.

(1)
Since $i^!$ and $R^1i^!$ vanish on all $Z$-modules,
(\ref{eq.five.seq}) gives $\Ext^1_X(i_*M,\b_*N)=0$.
This is precisely the criterion for $Z$ to be contained
in $X \backslash Y$.

(2)
Since $i^!$ vanishes on all $Z$-modules, the only $Z$-module
belonging to $\Mod_Y X$ is the zero module. Hence
$\Hom_X(\b_*N,-)$ is zero on $\Mod_Y X$.
By Lemma \ref{lem.VV}, all the terms in the minimal injective
resolution of $i_*M$ belong to $\Mod_Y X$; hence
$\Ext^n_X(\b_*N,i_*M)=0$ for all $n$. Hence $Y$ is contained in $X
\backslash Z$.
\end{pf}

\begin{corollary}
\label{cor.pts.in.divisor}
Let $Y$ be an effective divisor on $X$. A closed
point in $X$ is contained in either $Y$ or $X\backslash Y$ but
not both, and is a closed point in whichever of these spaces
contains it.
\end{corollary}
\begin{pf}
Let $p$ be a closed point of $X$ and suppose that $p$ does not lie
on $Y$. By Proposition \ref{prop.contained.in.divisor}, $p$ 
is contained in $X \backslash Y$.
As remarked in the proof of Corollary \ref{cor.pts.in.complement},
$p$ is closed in $X \backslash Y$ and hence is a closed point of $X
\backslash Y$.
\end{pf}

It is not true, even in the commutative case, that the set of
closed points in $X$ is the disjoint union of the set of closed
points in $Y$ and the set of closed points in $X \backslash Y$;
consider the spectrum of a discrete valuation ring. 
The next result gives conditions under which that is true.

\begin{proposition}
\label{prop.big.pt}
Let $X$ be a locally noetherian space. 
Suppose that $j:X \backslash Y \to X$ is the
inclusion of the complement to an effective divisor.
Let $p \in X \backslash Y$ be a closed $k$-rational point of this
complement.
If $j_*\cO_p$ is a noetherian $X$-module, then it is 
simple, and $j(p)$ (see Definition \ref{defn.pt.image})
is a closed $k$-rational point of $X$.
\end{proposition}
\begin{pf}
We will show first that if $j_*\cO_p$ is not a simple $X$-module,
then it is not noetherian. To that end, suppose that $M$ is a
non-zero $X$-submodule of $j_*\cO_p$ that is not equal to
$j_*\cO_p$. 

Let $i:Y \to X$ be the inclusion and let
$\tau:\Mod X \to \Mod_Y X$ be the $Y$-torsion functor.
Since $i^!j_*=0$, $i^!M=0$ and $\tau M=0$. Hence there is an exact
sequence
\begin{equation}
\label{eq.big.pt}
0 \to M \to j_*\cO_p \to R^1\tau M \to 0.
\end{equation}
Since $M$ is torsion-free, $j^*M \ne 0$, whence $j^*M \cong
j^*j_*\cO_p \cong \cO_p$. 
Because $i^*j_*=0$, applying (\ref{eq.V.crit}) to $j_*\cO_p$ gives
$j_*\cO_p \cong (j_*\cO_p)(Y)$. 
Then applying (\ref{eq.V.crit}) to the other
two terms in (\ref{eq.big.pt}) gives the other 
columns in the commutative diagram below:
$$
\begin{CD}
@. 0 @. 0 @. i^!(R^1\tau M)
\\
@. @VVV @VVV @VVV
\\
0 @>>> M @>>> j_*\cO_p @>>> R^1\tau M @>>> 0
\\
@. @VVV @VV{\cong}V @VVV
\\
0 @>>> M(Y) @>>> j_*\cO_p(Y) @>>> (R^1\tau M)(Y) @>>> 0
\\
@. @VVV @VVV @VVV
\\
@. (i^*M)(Y) @. 0 @. 0
\\
@. @VVV
\\
@. 0
\end{CD}
$$
Let $M_1$ denote the image of $M(Y)$ in $j_*\cO_p$.
Repeating this argument produces a chain of submodules
$$
M=M_0 \subset M_1 \subset M_2 \subset \cdots \subset j_*\cO_p
$$
such that $M_n \cong M(nY)$ and $M_n/M_{n-1} \cong (i^*M)(nY)$.
By the Snake Lemma, $(i^*M)(Y) \cong i^!(R^1\tau M)$; since
$R^1\tau M$ is non-zero, $i^!(R^1\tau M)
\ne 0$, whence $i^*M \ne 0$. Hence the chain of submodules is
strictly increasing, thus showing that $j_*\cO_p$ is not noetherian.

When $j_*\cO_p$ is noetherian there can be no such $M$, 
so $j_*\cO_p$ is simple.
Now we show that it is a tiny simple. First,
$\Hom_X(j_*\cO_p,j_*\cO_p) \cong \Hom_{X \backslash
Y}(\cO_p,\cO_p) \cong k$. Let $N$ be a noetherian
$X$-module. Then $j^*N$ is a noetherian $X\backslash Y$-module, 
so $\Hom_X(N,j_*\cO_p) \cong \Hom_{X \backslash Y}(j^*N,\cO_p)$ is
a finite dimensional $k$-vector space. Thus $j_*\cO_p$ is a tiny
simple, so there is a closed $k$-rational point $p' \in X$ such
that $j_*\cO_p =\cO_{p'}$.
\end{pf}

{\bf Remark.}
Even in the commutative case the hypothesis in Proposition
\ref{prop.big.pt} does not always hold: for example, let $X=\Spec
R$ where $R$ is a discrete valuation ring. 
However, if $X$ is a scheme of finite type over a
field $k$, then every closed $k$-rational point in the complement
of an effective divisor is a closed $k$-rational point of $X$. 
In the non-commutative case this sort of finiteness hypothesis 
does not ensure that a $k$-rational closed point of 
$X \backslash Y$ is a closed point of $X$.

For example, consider the graded line defined by $\Mod
\LL^1=\GrMod k[x]$ where $\deg x=1$. Then $x$ cuts out an effective
divisor $Y$ defined by $\Mod Y=\GrMod k[x]/(x)$, and $\LL^1 \backslash
Y$ is isomorphic to $\Spec k$ ($\Mod \LL^1 \backslash Y=\Tails k[x] =
\PP^0_k$). 
However, if $j:\LL^1 \backslash Y \to \LL^1$ denotes the inclusion, then
$j_*\cO_{\LL^1\backslash Y} \cong k[x,x^{-1}]$. Since $\LL^1 \backslash
Y$ is the generic point of $\LL^1$ (see \cite{SmInt}), it is no
surprise that it is not a closed point of $\LL^1$.

\begin{corollary}
\label{cor.big.pt}
Let $X$ be a locally noetherian space. 
Suppose that $j:X \backslash Y \to X$ is the
inclusion of the complement to an effective divisor.
If $j_*\cO_p$ is a noetherian $X$-module for every closed point $p
\in X \backslash Y$, then the set of closed points on $X$ is 
the disjoint union of the sets of closed
points on $Y$ and on $X\backslash Y$.
\end{corollary}
\begin{pf}
It is clear that a closed point on $Y$ is a closed point on $X$, so
the result follows from Proposition \ref{prop.big.pt} and Corollary
\ref{cor.pts.in.divisor}.
\end{pf}

Another way of stating the conclusion of this is that $j_*$ and
$j^*$ set up a bijection between the closed points of $X$ that do
not lie on $Y$, and the closed points of $X \backslash Y$.

\medskip

Effective divisors on a space typically arise from regular elements in 
a (homogeneous) coordinate ring. For example, if $z$ is a regular
element in a ring $R$ such that $zR=Rz$, then $\Mod R/(z)$ is an
effective divisor on $\Mod R$.

\begin{proposition}
\label{prop.divisor}
Let $A$ be a graded ring with a central regular element $z \in A_d$.
Set $X=\Tails A$, $H=\Tails A/(z)$, and let $i:H \to X$ be the
inclusion. 
Then $H$ is an effective divisor on $X$ and $M(H) \cong M(d)$
for all $M \in \Mod X$.
Furthermore, if $A$ is connected and generated in degree one, then
$X \backslash H$ is the affine space with coordinate ring
$A[z^{-1}]_0$.
\end{proposition}
\begin{pf}
We use the notation and terminology in \cite{vdB}.

Define the bimodule $o_H$ by $\cHom_X(o_H,-)=i_*i^!$. 
Thus $-\otimes o_H \cong i_*i^*.$ 
As usual $o_X$ denotes the identity functor on
$\Mod X$ considered as a bimodule. Define $o_X(-H)$ as the kernel 
of the map of bimodules $o_X \to o_H$ that corresponds to the counit
$i_*i^! \to \id_X$. Thus $\cHom_X(o_X(-H),-)$ is defined by the requirement
that if $I$ is an injective $X$-module, then  $\cHom_X(o_X(-H),I)$ 
is the cokernel in the exact sequence
$$
0 \to \cHom_X(o_H,I) \to \cHom_X(o_X,I) \to \cHom_X(o_X(-H),I) \to 0.
$$
This sequence can also be written as 
\begin{equation}
\label{eq.I.seq}
0 \to i_*i^! I \to I \to \cHom_X(o_X(-H),I) \to 0.
\end{equation}
Let $\pi:\GrMod A \to \Mod X$ be the quotient functor.
There is an injective $Q$ in $\GrMod A$ such that $I=\pi Q$. 
It is easy to
see that there is an exact sequence 
$$
0 \to K \to Q @>{\c}>> Q(d) \to 0
$$
where $\c(q)=qz$. Since $\pi K \cong i_*i^!I$, applying $\pi$ to this 
gives (\ref{eq.I.seq}). Therefore $\cHom_X(o_X(-H),I) 
\cong I(d)$. This shows that $\cHom_X(o_X(-H),-)$ is an invertible
functor, and hence that $H$ is an effective divisor on $X$. The bimodule
$o_X(H)$ is defined by the requirement that $\cHom_X(o_X(H),-)$ is an
inverse to $\cHom_X(o_X(-H),-)$, so it follows that $I \otimes o_X(H)
\cong I(d)$. Now taking an injective resolution of an arbitrary $X$-module,
we see that $M(H) \cong M(d)$ for all $M \in \Mod X$.

By definition $\Mod X \backslash H =\Mod X/\Mod_HX \cong \GrMod
A/\sT$ where $\sT$ is the localizing subcategory consisting of the
$A$-modules which are $z$-torsion.
Thus $\GrMod A/\sT \cong \GrMod A[z^{-1}]$. When $A$ is
connected and generated by $A_1$, $A[z^{-1}]$ is strongly graded,
meaning that $1$ belongs to the product of the degree $n$ and
degree $-n$ components for all $n$; this implies that $ \GrMod
A[z^{-1}]$ is equivalent to $\Mod A[z^{-1}]$ via the functor that
takes the degree zero component of a module.
\end{pf}

An affine scheme can be embedded in a projective scheme
in such a way that it is the open complement of an ample
divisor. It is well-known that a non-commutative affine 
space can be embedded in a non-commutative projective space
in a similar way. We show next that closed points behave well with
respect to this embedding.

First we recall the embedding just alluded to.
Let $R$ be an algebra over a field $k$. Fix a set of
generators for $R$ and write $R_n$ for the $k$-linear span of all 
products of $\le n$ of those generators. 
Let $t$ denote an indeterminate commuting with $R$, and set
$$
\tilde R:=R_0 \oplus R_1t \oplus R_2 t^2 \oplus \cdots \subset
R[t].
$$
Define
$$
U=\Mod R 
\qquad
\hbox{and} 
\qquad
X=\Tails \tilde R.
$$ 
We write 
$$
\pi:\GrMod \tilde R \to \Mod X 
\qquad \hbox{and} \qquad 
\omega: \Mod X \to \GrMod \tilde R
$$ 
for the quotient functor and its right adjoint.

The closed subspace of $X$ where $t\in \tilde R_1$ vanishes is
denoted by $H$ and is called the {\sf hyperplane at infinity}.
Proposition \ref{prop.divisor} shows that $H$ is an effective divisor
because $\Mod H \cong \Mod \tilde R/(t).$ We also have 
$\Mod H \cong \Tails \gr R$, 
where $\gr R$ is the associated graded ring with respect to the 
filtration $R_n$. 
Because $R[t^{-1}]_0 \cong R$, there is a natural map 
$$
j:U \to X
$$
sending $U$ isomorphically to $X \backslash H$.

\begin{proposition}
\label{prop.Rees.pts}
Suppose that $\tilde R$ is a finitely generated right noetherian 
$k$-algebra.
The functors $j_*$ and $j^*$ set up a bijection between
\begin{enumerate}
\item{}
the finite dimensional simple $R$-modules $V$ such 
that $\End_R V \cong k$ and 
\item{}
the closed $k$-rational points of $X$ that do not lie on $H$.
\end{enumerate}
\end{proposition}
\begin{pf}
This will follow from Corollary \ref{cor.big.pt} once we show that
$j_*V$ is a noetherian $X$-module for every finite dimensional
simple right $R$-module $V$.
By definition, $j_*V \cong \pi \tilde V$, where $\tilde V$ is
the graded $\tilde R$-module
$$
V \oplus Vt \oplus Vt^2 \oplus \ldots
$$
with degree $n$ component $\tilde V_n =Vt^n$ and action
$$
(vt^n).(rt^i):=(vr)t^{n+i}
$$
for $v \in V$ and $r \in R_i$.
Since every component of $\tilde V$ has the same finite dimension,
$\pi \tilde V$ is noetherian.
\end{pf}

\begin{lemma}
\label{lem.Ext1=0}
Let $i:Y \to X$ be the inclusion of effective divisor in a space $X$.
Let $j:X \backslash Y \to X$ be the inclusion of the weakly open complement.
Then $\Ext^1_X(j_*N,M)=0$ for all $N \in \Mod X \backslash Y$ and 
$M \in \Mod Y$.
\end{lemma}
\begin{pf}
If we apply $i^*$ to an exact sequence
\begin{equation}
0 \to M @>{\a}>> K @>>> j_*N \to 0
\end{equation}
we obtain an exact sequence
$$
(L_1i^*)(j_*N) @>>> i^*M @>{i^*(\a)}>> i^*K @>>> i^*j_*N \to 0.
$$
The last term of this is zero by Proposition \ref{prop.open}, and the
first is zero because $i^!j_*=0$. Hence there is a map
$\b:i^*K \to i^*M$ such that $\b \circ i^*(\a)=\id_M$. However, 
$i^*(\a)=\c \circ \a$ where $\c:M \to i^*M$ is the canonical epimorphism,
so $\b\c\a=\id_M$, which shows that the original exact sequence splits.
Hence the result.
\end{pf}

\medskip

The next result is an analogue of the result that an ample divisor
on a scheme meets every subscheme of positive dimension.

Let $X$ be a locally noetherian space and $\mod X$ its category of
noetherian modules. We denote the Krull
dimension, in the sense of Gabriel \cite{G}, of $M \in \mod X$
by $\Kdim M$. We say that a noetherian $X$-module $M$ is
$d$-critical if $\Kdim M=d$ but $\Kdim M/N < d$ for all $0 \ne N
\subset M$. A noetherian module of Krull dimension $d$ has a
$d$-critical quotient.
We define the {\sf dimension} of a weakly closed subspace $W$ of
$X$ to be 
$$
\dim W:=\max\{ \Kdim M \; | \; M \in \mod W\}.
$$

\begin{definition}
\label{defn.ample}
Let $X$ be a locally noetherian space endowed with a structure map
$\pi:X \to \Spec k$. Set $\cO_X=\pi^*k$ and
$H^q(X,-)=\Ext^q_X(\cO_X,-)$.
An effective divisor $Y$ on $X$ is
{\sf ample} if for every $M \in \mod X$
\begin{itemize}
\item{}
there are integers $n_i \le 0$ and an
epimorphism $\bigoplus \cO_X(n_iY) \to M$, and
\item{}
$H^q(X,M(nY))=0$ if $q \ne 0$ and $n \gg 0$.
\end{itemize}
\end{definition}

\begin{proposition}
\label{prop.ample.nef}
Let $X$ be a locally noetherian space over $\Spec k$, and let $Y$
be an ample effective divisor on $X$. Suppose that 
$\dim_k \Ext_X^n(M,N)< \infty$ for all $M, N \in \mod X$. 
If $W$ is a weakly closed subspace of $X$ and $\dim W>0$,
then $W \cap Y \ne \phi$.
\end{proposition}
\begin{pf}
We will assume that $W \cap Y=\phi$ and get a contradiction.
Let $d=\dim W$ and let $M$ be a $d$-critical noetherian $W$-module.
By (\ref{eq.V.crit}), there is an exact
sequence $0 \to i^!M(-Y) \to M(-Y) \to M \to i^*M \to 0$.
Since $W \cap Y=\phi$, $i^!M=i^*M=0$, whence $M \cong M(Y)$. 
Hence $\dim_k H^0(X,M(nY))$ is a finite number independent of $n$. 
By \cite{AZ}, $\Mod X$ is equivalent to $\Tails B$ where $B$ is the 
graded algebra $\oplus_{n=0}^\infty H^0(X,\cO_X(nY))$. Under this
equivalence, $M \cong \pi V$, where $V$ is the 
graded $B$-module 
$\oplus_{n=0}^\infty H^0(X,M(nY))$. But $\dim V_n$ is
finite and independent of $n$ so $\pi V$ has finite length. 
This contradicts our hypothesis that $\Kdim M =d \ge 1$.
\end{pf}

The next result generalizes the surprising result in 
\cite[Lemma 6.4]{KKO}. In that
result, $Y$ is an effective divisor and $X \backslash Y$ is an 
affine space having a coordinate ring that is an infinite
dimensional $\CC$-algebra; such an algebra cannot have a finite
dimensional simple module, so $X \backslash Y$ has no closed
points.

\begin{proposition}
Let $k$ be an algebraically closed field and $X$ a locally
noetherian $\Ext$-finite space over $k$. Let $i:Y \to X$ be the
inclusion of an effective divisor. Suppose that $X \backslash Y$ 
has no closed $k$-rational points. If $\cF$ is a noetherian
$X$-module such that $i^*\cF=0$, then $\cF=0$.
\end{proposition}
\begin{pf}
Suppose that $\cF \ne 0$.
Since $\cF$ is noetherian it has a simple quotient. By 
Corollary \ref{cor.one.pt}, that simple is of the form $\cO_p$ for
a closed $k$-rational point $p \in X$. By the hypothesis and
Corollary \ref{cor.pts.in.divisor}, $p \in Y$; this implies that
$i^*\cF \ne 0$.
\end{pf}


\begin{thebibliography}{11}


\bibitem{ATV1} M. Artin, J. Tate and M. van den Bergh, Some
algebras related to automorphisms of elliptic curves,
in {\it The
Grothendieck Festschrift, Vol.1,} 33-85, Birkhauser, Boston 1990.

\bibitem{ATV2}
M. Artin, J. Tate and M. van den Bergh,
Modules over regular algebras of dimension 3, {\it Invent. Math.,}
{\bf 106} (1991) 335-388.


\bibitem{AZ}
M. Artin and J. Zhang,
Non-commutative Projective Schemes,
{\it Adv.\ Math.}, {\bf 109} (1994), 228--287.

\bibitem{C}
G. Cauchon,
Les T-anneaux, la condition (H) de Gabriel et ses cons\'equences,
{\it Comm. Alg.,} {\bf 4} (1976) 11-50.

\bibitem{G}

P.\ Gabriel, Des cat\'{e}gories ab\'{e}liennes, {\it Bull.\ Soc.\
  Math.\ France} {\bf 90} (1962), 328--448.

\bibitem{G.SGA6}
A. Grothendieck (with P. Berthelot and L. Illusie), 
{\it Th\'eorie des intersections et Th\'eor\`eme de
Riemann-Roch} (SGA 6), 
Lecture Notes in Math. 225, Springer-Verlag, Berlin/New York, 1971.

\bibitem{H.RD}
R. Hartshorne, {\it Residues and Duality,}
Lecture Notes in Math. 20, Springer-Verlag, Berlin/New
York, 1966.


\bibitem{H}
R. Hartshorne, {\it Algebraic Geometry,} Springer-Verlag, 1977.





\bibitem{J}
P. J\o rgensen,
Intersection theory on non-commutative surfaces, 
{\sl Trans. Amer.  Math. Soc.,} {\bf 352} (2000) 5817-5854.

\bibitem{JS1}
P. J\o rgensen and S.P. Smith,
Curves on non-commutative surfaces, in preparation.

\bibitem{JS2}
P. J\o rgensen and S.P. Smith,
Intersection theory for the blow-up of a non-commutative
surface, in preparation.

\bibitem{JS}
A. Joyal and R. Street, An introduction to Tanaka duality and
quantum groups, pp. 413-492 in 
{\it Category Theory, Proceedings, Como 1990,}
Eds. A. carboni, M.C. Pedicchio, G. Rosolini,
Lect. Notes in Math., vol. 1488, Springer Verlag, 1991.

\bibitem{KKO}
A. Kapustin, A. Kuznetzov, and D. Orlov,
Noncommutative Instantons and Twistor Transform,
Preprint (2000) hep-th/0002193 v3.

\bibitem{MR}
J.\ C.\ McConnell and J.\ C.\ Robson,
{\it Noncommutative Noetherian Rings,}
John Wiley and Sons, 1987.

\bibitem{Papp}
C. Pappacena, 
The injective spectrum of a non-commutative space,
manuscript, February 2001.

\bibitem{Pop}
N. Popescu,
{\it Abelian Categories with Applications to Rings and Modules,}
Academic Press (Lond. Math. Soc. Monographs), 1973.



\bibitem{Rosen}
A.L. Rosenberg, {\it Non-commutative algebraic geometry and
representations of quantized algebras,} Vol. 330, Kluwer
Academic Publishers, 1995.

\bibitem{Ros2}
A.L. Rosenberg, Non-commutative schemes, {\it Compos. Math.,}
{\bf 112} (1998) 93-125.

\bibitem{SZ}
S.P. Smith and J. Zhang, Curves on quasi-schemes, 
{\it Algebras and Representation Theory,}
{\bf 1} (1998) 311-351.


\bibitem{SmInt}
S.P. Smith,
Integral non-commutative spaces,
{\it J. Algebra,} to appear.

\bibitem{SmCo}
S.P. Smith,
Comodules and non-commutative spaces,
in preparation.


\bibitem{SmMaps}
S.P. Smith,
Maps between non-commutative spaces,
Preprint, University of Washington, 2000.



\bibitem{vdB}
M. Van den Bergh, Blowing up on non-commutative smooth surfaces,
{\it Mem. Amer. Math. Soc.,} to appear.

\bibitem{VV}
M.\ Van Gastel and M.\ Van den Bergh,
Graded modules of Gelfand-Kirillov dimension
one over three-dimensional Artin-Schelter regular algebras,
{\it J.\ Algebra}, {\bf 196} (1997), 697--707.


\bibitem{Ver}
A.B. Verevkin, On a non-commutative analogue of the
category of coherent sheaves on a projective scheme, {\em Amer. Math.
Soc.  Trans.,} {\bf 151} (1992) 41-53.


\end{thebibliography}
\end{document}